\documentclass[lualatex,sn-mathphys-num]{sn-jnl}

\usepackage{graphicx}%
\usepackage{multirow}%
\usepackage{amsmath,amssymb,amsfonts}%
\usepackage{amsthm}%
\usepackage{mathrsfs}%
\usepackage[title]{appendix}%
\usepackage{xcolor}%
\usepackage{textcomp}%
\usepackage{manyfoot}%
\usepackage{booktabs}%
\usepackage{algorithm}%
\usepackage{algorithmicx}%
\usepackage{algpseudocode}%
\usepackage{listings}%
\usepackage{fullpage}%

\newcommand{\bcell}{$\beta$-cell}
\newcommand{\bcells}{$\beta$-cells}
\newcommand{\Ca}{Ca$^{2+}$}
\newcommand{\Kp}{K$^+$}

\usepackage{comment}

\usepackage{mathtools}

\usepackage{tabularray}
\SetTblrInner[tblr,talltblr]{rowsep=0.75ex} 

\usepackage{siunitx}
\UseTblrLibrary{booktabs}
\UseTblrLibrary{siunitx}

\usepackage{amsopn}

\DeclareMathOperator{\tr}{tr}
\DeclareMathOperator{\adj}{adj}
\DeclareMathOperator{\rk}{rk}
\DeclareMathOperator{\dist}{dist}

\usepackage{enumitem}

\usepackage[capitalise, noabbrev]{cleveref}
\hypersetup{hypertexnames=false}
\crefname{appendix}{Appendix}{Appendices}
\Crefname{appendix}{Appendix}{Appendices}

\theoremstyle{thmstyleone}%
\newtheorem{theorem}{Theorem}
\newtheorem{proposition}[theorem]{Proposition}%

\theoremstyle{thmstyletwo}%
\newtheorem{example}{Example}%
\newtheorem{remark}{Remark}%

\theoremstyle{thmstylethree}%
\newtheorem{definition}{Definition}%

\begin{document}

\title[]{Mathematical analysis of a model for glycolytic oscillations in pancreatic $\beta$-cells}


\author*[1]{\fnm{Prannath} \sur{Moolchand}}\email{prannath.moolchand@sydney.edu.au}

\author[2, 3]{\fnm{Richard} \sur{Bertram}}\email{rbertram@fsu.edu}

\author[1]{\fnm{Martin} \sur{Wechselberger}}\email{martin.wechselberger@sydney.edu.au}

\affil*[1]{\orgdiv{School of Mathematics and Statistics}, \orgname{The University of Sydney}, \orgaddress{\street{Eastern Avenue}, \city{Camperdown}, \postcode{2006}, \state{NSW}, \country{Australia}}}

\affil[2]{\orgdiv{Department of Mathematics}, \orgname{Florida State University}, \orgaddress{\street{1017 Academic Way}, \city{Tallahassee}, \postcode{32306-4510}, \state{Florida}, \country{United States of America}}}

\affil[3]{\orgdiv{Institute of Molecular Biophysics}, \orgname{Florida State University}, \orgaddress{\street{91 Chieftan Way}, \city{Tallahassee}, \postcode{32306-4380}, \state{Florida}, \country{United States of America}}}


\abstract{
Insulin secretion is pulsatile, with a period of about 5 min, coexisting with an ultradian rhythm (period of about 1 h), and a circadian (daily) rhythm. These rhythms are disrupted in type II diabetes. The 5-min rhythm is typically maintained even at fasting glucose levels, when there is no electrical activity in the insulin-secreting pancreatic $\beta$-cells. It has been proposed that this rhythm is due to oscillations in glycolysis, the first stage of glucose metabolism. A mathematical model for these oscillations was developed by Paul Smolen (the `Smolen model') and is an important element of the Integrated Oscillator Model (IOM) for $\beta$-cell activity across all levels of glucose.
In this article, we provide a detailed mathematical analysis of the dynamics underlying oscillations in the Smolen model. In particular, we demonstrate how the biophysical mechanism for oscillations in this model, a substrate-depletion mechanism, can be decomposed mathematically using geometric singular perturbation analysis. We link the `pulse' and `refill' stages of the oscillation to distinct dynamic regimes and describe the flow dynamics through these regimes. This is the first detailed mathematical analysis of the Smolen model for glycolytic oscillations that are hypothesised to occur at fasting glucose levels and, in some cases, under superthreshold conditions that would occur after a meal.  
}




\maketitle

\section{Introduction}
Insulin is a hormone promoting glucose uptake into cells, and its secretion by pancreatic $\beta$-cells of the islets of Langerhans is fundamental to metabolic health, with dysfunction directly linked to diabetes \cite{ashcroft13,satin15}. $\beta$-cells share key electrophysiological traits with neurons, exhibiting electrical excitability due to a wide variety of ion channels in the plasma membrane \cite{bertram07b}. Thus, the modelling of $\beta$-cell dynamics has a rich history, evolving from the Chay-Keizer model \cite{Chay1983}, which itself descended from the biophysically grounded Hodgkin-Huxley model \cite{Hodgkin1952}, to the Integrated Oscillator Model (IOM), which contains modules for metabolism, \Ca\ handling, and electrical activity \cite{Bertram2023}.

A distinguishing feature of \bcells\ compared to the majority of other electrically excitable cells is the presence of ion channels regulated by the cell's metabolic state. These ATP-sensitive \Kp\ channels are inactivated by intracellular adenosine triphosphate (ATP), which is the end product of glucose metabolism \cite{ashcroft13}. The ionic current produced through these channels links the cell's metabolic state to its electrical activity. This electrical activity brings \Ca\ into the cell, which causes the exocytosis of insulin. Electrical activity comes in a bursting pattern, so insulin is released from the cells into the bloodstream in pulses \cite{ashcroft13,satin15}.

\begin{figure}[ht]
\centering
\includegraphics[width=0.6\textwidth]{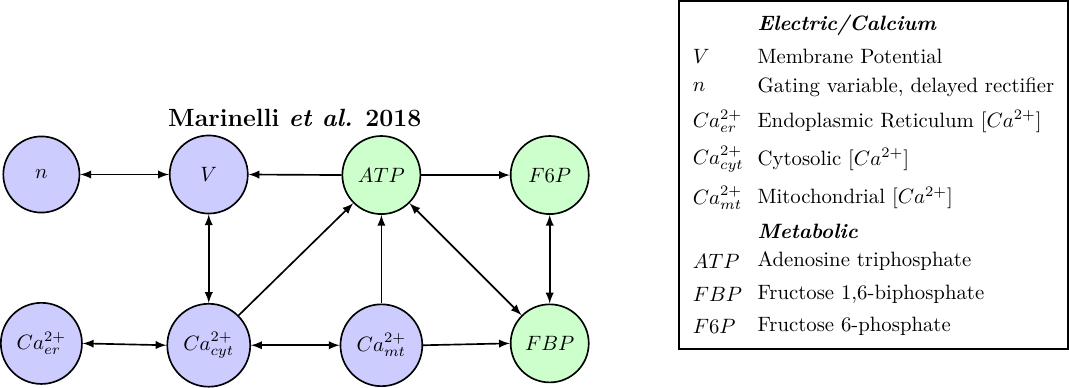}
\caption{Network diagram of the 8-dimensional Integrated Oscillator Model (IOM) \cite{Marinelli2018}, showing the interaction between electrical/\Ca\ (blue) and metabolic (green) subsystems.}
\label{fig:coupling-sketch}
\end{figure}

The bursting electrical activity occurs when the glucose concentration is about 6 mM or higher, which is the case after a meal. However, between meals, the fasting glucose level is lower and there is little or no electrical activity. Yet, rhythmic insulin levels are often observed even at these fasting glucose levels \cite{goodner77,lang_d79,porksen96a,song_sh00}. It has been proposed that the rhythmic insulin secretion at these fasting glucose levels is due to oscillations in glycolysis \cite{fletcher22}. These oscillations are produced through the action of the allosteric enzyme phosphofructokinase (PFK) \cite{tornheim75}. The muscle isoform of PFK dominates glycolytic activity in \bcells\ \cite{yaney95}, and it reliably generates oscillations. The prospect that glycolytic oscillations may control oscillations in pancreatic \bcells\ led Paul Smolen to develop a mathematical model for glycolytic oscillations generated by this isoform of the PFK enzyme \cite{Smolen1995}. This model was later incorporated into the Integrated Oscillator Model (IOM) that combines \bcell\ metabolism, intracellular \Ca\ handling, and electrical activity (\cref{fig:coupling-sketch}) \cite{Bertram2023,bertram07a,Marinelli2018}. This is the most comprehensive model for oscillations in \bcells, and the glycolytic component is key to understanding oscillations at fasting glucose levels as well as compound bursting oscillations (episodes of bursts) that are often seen at superthreshold glucose levels \cite{bertram04c}. 

Glycolytic oscillations produced by the Smolen model can be described biophysically as arising from a substrate-depletion mechanism. The substrate of the enzyme, fructose-6-phosphate (F6P), produces the product fructose-1,6-bisphosphate (FBP), which can bind to an allosteric site on the PFK enzyme, amplifying production of the product. That is, the FBP binding to PFK produces positive feedback. The result is a pulse in FBP levels at the expense of the substrate F6P. During the pulse, F6P levels decline so that eventually there is not enough to sustain further FBP production. At this point, the FBP level crashes, while the F6P level slowly rises. Once it has risen far enough, a new cycle begins.

This pulse--refill cycle exhibits the fundamental structure of a \emph{two-stroke} oscillator. Following Le Corbeiller \cite{LeCorbeiller1960}, continuous oscillators can be classified by their energetic exchanges independently of timescale separation. \emph{Four-stroke} systems receive energy twice per period via a characteristic with two distinct bends, whereas \emph{two-stroke} systems receive energy only once per period via an asymmetric or single-bend characteristic. When such systems evolve on multiple timescales, standard models like the van der Pol \cite{vdP1926} and FitzHugh \cite{FitzHugh1961} equations manifest as four-stroke relaxation oscillators, comprising two slow and two fast phases. While the native Smolen model is energetically two-stroke, it lacks explicit fast-slow scale separation. However, applying a state-dependent time transformation reveals its formal multi-scale structure. The system then maps perfectly to the non-standard singular perturbation problems formalized by Jelbart and Wechselberger \cite{Jelbart2020} as \emph{two-stroke relaxation oscillators}---a class requiring only a single turning point in the characteristic to produce exactly one slow and one fast phase.

The central aim of this manuscript is to develop a rigorous mathematical understanding of these two-stroke glycolytic oscillations. By leveraging geometric singular perturbation theory (GSPT) \cite{Bertram2017}, we systematically unpack the multi-scale dynamics of the transformed Smolen model.

The outline of the paper is as follows: in \cref{sec:model}, we describe the Smolen model and perform non-dimensionalisation, identifying the key small parameter $\varepsilon$. We then apply a coordinate change and time-rescaling to transform the system into a form
that is more suitable for GSPT analysis. In \cref{sec:ana}, we conduct a GSPT analysis by defining three distinct scaling regimes to heuristically construct the singular limit cycle. 
In \cref{sec:blowup}, we formally prove the existence of the limit cycle by desingularising the flow via two successive geometric blow-ups. We conclude in \cref{sec:discussion} with a discussion of our findings, in particular, relating the different dynamic regimes to the pulse and refill stages of the Smolen model.

\section{The Smolen model}
\label{sec:model}
The Smolen model for glycolytic oscillations consists of three ordinary differential equations, for the variables F6P (the PFK substrate), FBP (the PFK product) and ATP (a PFK modulator). The first two equations, for F6P and FBP, are central to the production of oscillations. The third variable, ATP, is an inhibitor of the PFK enzyme, but is not part of the oscillation mechanism. In applications, the right-hand side of the ATP differential equation would reflect glycolytic and mitochondrial ATP production and various sources of ATP hydrolysis that convert ATP to ADP. Since we are studying the glycolytic oscillator in isolation from other cell components, and since ATP is not part of the oscillation mechanism, we treat it as a fixed parameter, as was done previously \cite{Bertram2023,bertram04c}. The model differential equations, based on \cite{Marinelli2018}, are then:
\begin{equation}
\begin{rcases}
\dfrac{d\text{F6P}}{dt} &= \dfrac{1}{1+K_{\rm GPI}} (J_{\rm GK} - J_{\rm PFK}) \\[1ex]
\dfrac{d\text{FBP}}{dt} &= \dfrac{1}{1+K_{\rm LG}} (J_{\rm PFK} -\frac{1}{2}J_{\rm PDH})
\end{rcases} \label{app:eq:f6pfbp} \noeqref{app:eq:f6pfbp}
\end{equation}
where the units of F6P and FBP are $\mu$M.
$J_{\rm GK}$ is the glucokinase flux that causes F6P production. $J_{\rm PFK}$ is the flux through the PFK enzyme, and is the entry point of the ATP parameter. PFK converts F6P to FBP. Finally, $J_{\rm PDH}$ is the flux of downstream glycolytic elements through the pyruvate dehydrogenase reaction. $K_{\rm GPI}$ and $K_{LG}$ are parameters. Other parameters enter the equations through the flux terms, which are written out in \cref{app:sec:bio} and described in detail in \cite{Marinelli2018}.

We rewrite the F6P and FBP equations as:
\begin{align}
&\begin{rcases}
\dfrac{dx}{dt} &= \beta \left(\alpha - \nu r(x,y;K)\right)\\[1ex]
\dfrac{dy}{dt} &= \eta \left( \nu r(x,y;K) - \gamma \sqrt{\dfrac{y}{\omega}}\right)
\end{rcases} \label{eq:sysxy} \noeqref{eq:sysxy}\\
r(x,y;K) &= \dfrac{\dfrac{x^2y}{\kappa_5} + \dfrac{x^2}{\kappa_6}}{\dfrac{x^2y}{\kappa_1} + \dfrac{x^2}{\kappa_2} + \dfrac{y}{\kappa_3} + \dfrac{1}{\kappa_4}} \label{eq:brp} \noeqref{eq:brp}
\end{align}
with state variables $x=\;$F6P and $y=\;$FBP, given in units of $\mu$M, and system parameters $(\alpha, \beta, \eta, \nu,\omega)$ and $K=(\kappa_1,\ldots,\kappa_6)$. 

\begin{figure}[t]
\centering
\includegraphics[width=0.9\textwidth]{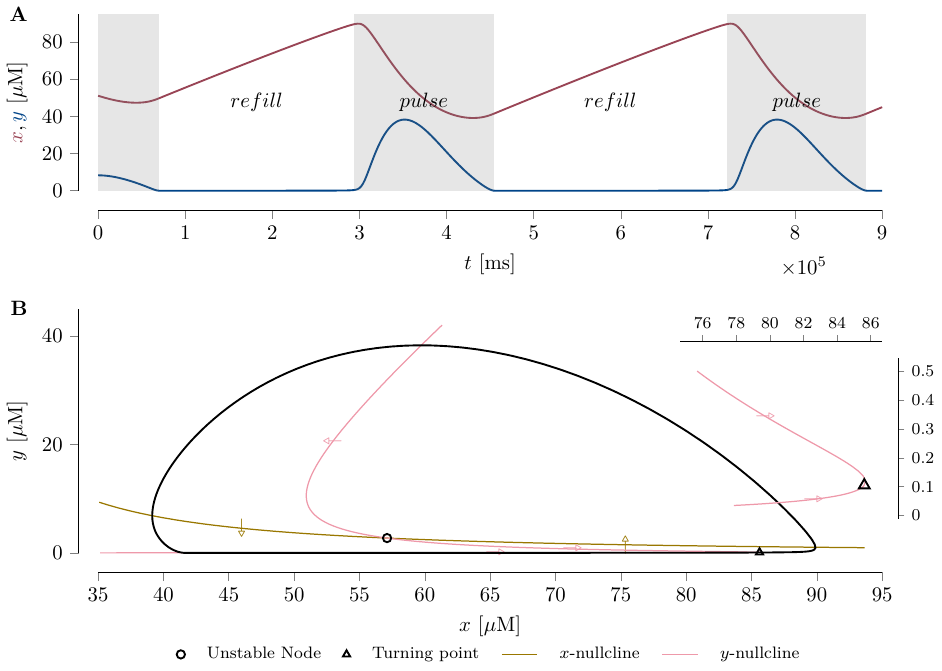}
\caption{Dynamics of oscillations in the Smolen model \eqref{eq:sysxy}. \textbf{A:} Time courses show a slow build-up of $x$ during the refill stage  followed by a pulse of $y$ during the pulse phase.  \textbf{B:} Phase portrait showing nullclines and the limit cycle for the glycolytic oscillations  (solid black curve).}
\label{fig:origprob}
\end{figure}

In \cref{fig:origprob}, the time courses of the glycolytic model \eqref{eq:sysxy} show the structure of a `two-stroke' oscillator \cite{Jelbart2020} comprising a `refill' phase when the enzyme substrate $x=$ F6P refills and a `pulse' phase when there is a large pulse in the enzyme product $y=$ FBP. In the $xy$ phase plane, the trajectory follows the bottom branch of the cubic-like $y$-nullcline during the refill phase, then moves past the right fold point and makes a large excursion during the pulse phase. The fold point or turning point of the $y$-nullcline is enlarged in the inset to make it more visible. The time courses and phase portrait both indicate that there are multiple timescales in the limit cycle dynamics. 

\subsection{Non-dimensionalisation and hierarchy of small parameters}
We next non-dimensionalise the equations and identify key parameters. This is done by rescaling the dependent and independent variables ($x,y,t$): 
\begin{equation}
x = \kappa_X X, \qquad
y = \kappa_Y Y, \qquad
t = \kappa_\tau\tau, \qquad
 \label{eq:sctr} \noeqref{eq:sctr}
\end{equation}
where $X$, $Y$ and $\tau$ are the new dimensionless quantities with $\kappa_X$, $\kappa_Y$ and $\kappa_\tau$ as reference scales with the same units as their biophysical counterparts. This gives the corresponding dimensionless equations for the Smolen model:

\begin{align}
&\begin{rcases}
\dfrac{dX}{d\tau} &= \hat{\alpha} - \hat{\nu}_1 \hat{r}(X,Y;\Sigma)\\[1ex]
\dfrac{dY}{d\tau} &= \hat{\nu}_2 \hat{r}(X,Y;\Sigma) - \hat{\gamma} \sqrt{Y}
\end{rcases} \label{eq:rndxy} \noeqref{eq:rndxy}\\
\hat{r}(X,Y; \Sigma) & := \frac{X^2 Y + \hat\sigma_6 X^2}{\hat\sigma_1 X^2 Y + \hat\sigma_2 X^2 + \hat\sigma_3 Y + \hat\sigma_4} \label{eq:rhatsig} \noeqref{eq:rhatsig}
\end{align}
with dimensionless rational function $\hat{r}(X,Y; \Sigma)$, and dimensionless parameters $\Sigma=(\hat\sigma_1,\ldots,\hat\sigma_4,\hat\sigma_6)$ and $(\hat\alpha,\hat\nu_1,\hat\nu_2,\hat\gamma)$; see \cref{tab:adimlspar} for the algebraic expressions of all dimensionless parameters.\\ 

\begin{talltblr}[
  caption = {Dimensionless Parameters for System \eqref{eq:rndxy}.}, %
  label   = {tab:adimlspar},                                         %
]{
  colspec = {c} 
}
$\begin{tblr}{
  colspec = {c|c||c|c||c|c||c|c||c|c},
  columns = {mode=math} %
}
\text{Par.} & \text{Expr.} & \text{Par.} & \text{Expr.} & \text{Par.} & \text{Expr.} & \text{Par.} & \text{Expr.} & \text{Par.} & \text{Expr.}\\ 
\hline
\hat\sigma_1 & \dfrac{\kappa_5}{\kappa_1}                           &  
\hat\sigma_2 & \dfrac{\kappa_5}{\kappa_2 \kappa_Y}                    &        
\hat\sigma_3 & \dfrac{\kappa_5}{\kappa_3 \kappa_X^2}                  & 
\hat\sigma_4 & \dfrac{\kappa_5}{\kappa_4 \kappa_X^2 \kappa_Y}         &      
\hat\sigma_6 & \dfrac{\kappa_5}{\kappa_6 \kappa_Y}                    \\  
\hat{\alpha}  & \kappa_\tau \dfrac{\beta \alpha}{\kappa_X}        &
\hat{\nu}_1   & \kappa_\tau \dfrac{\beta \nu}{\kappa_X}           & 
\hat{\nu}_2   & \kappa_\tau \dfrac{\eta \nu}{\kappa_Y}            & 
\hat{\gamma}  & \kappa_\tau \dfrac{\eta \gamma}{\sqrt{\kappa_Y \omega}} &
              & 
\end{tblr}$ \\
\end{talltblr}\\

To identify a key small parameter for a perturbation analysis, we selected the following scales:
\begin{equation}
\kappa_X =  \left( \dfrac{\kappa_5}{\kappa_4} \right)^\frac{1}{3}
\left( \dfrac{\kappa_6^3}{\kappa_4 \kappa_5^2} \right)^{\frac{1}{12}}, \qquad
\kappa_Y = \left( \dfrac{\kappa_5}{\kappa_4} \right)^\frac{1}{3}, \qquad
\kappa_\tau = \dfrac{1}{\eta \nu} \left( \dfrac{\kappa_5}{\kappa_4} \right)^\frac{1}{3}\,.
\label{eq:sclfac} \noeqref{eq:sclfac}
\end{equation}
This choice leads to the dimensionless parameter values summarised in \cref{tab:sumparams}, and it identifies $\hat\sigma_6= \num{2.52e-3} \ll 1$ as the key small parameter. \\
%

\begin{talltblr}[
  caption = {Scaling factors and parameters for non-dimensionalised System \eqref{eq:rndxy}.}, %
  label   = {tab:sumparams},                                         %
]{
  colspec = {c} 
}
$\begin{tblr}{
 colspec = {c|c|l||c|c|l},
 columns = {mode=math} %
}
\text{Par} & \text{Expr.} & \text{Value} & \text{Par} & \text{Expr.} & \text{Value}\\ 
\hline
\hat{\alpha}  & \dfrac{\beta \alpha}{\eta \nu}  \left( \dfrac{\kappa_4 \kappa_5^2}{\kappa_6^3} \right)^{\frac{1}{12}} & \num{5.56e-03}
& 
\hat\sigma_1 & \dfrac{\kappa_5}{\kappa_1}                                     & \num{1.38e+00}
\\
\hat{\nu}_1   & \dfrac{\beta}{\eta}  \left( \dfrac{\kappa_4 \kappa_5^2}{\kappa_6^3}\right)^{\frac{1}{12}} & \num{6.72e-02}
&
\hat\sigma_2 & \left(\dfrac{\kappa_6}{\kappa_2}\right) \left( \dfrac{\kappa_4 \kappa_5^2}{\kappa_6^3}\right)^{\frac{1}{3}}      & \num{4.20e-02}
\\
\hat{\nu}_2   & 1                                                                                                   & \num{1.00e-00}
&
\hat\sigma_3 & \left( \dfrac{\kappa_4^2 \kappa_5}{\kappa_3^3}\right)^{\frac{1}{3}} \left( \dfrac{\kappa_4 \kappa_5^2}{\kappa_6^3}\right)^{\frac{1}{6}} & \num{1.14e-01}
\\
\hat{\gamma}  & \dfrac{\gamma}{\nu \sqrt{\omega}} \left( \dfrac{\kappa_5}{\kappa_4} \right)^\frac{1}{6}               & \num{3.30e-01}
&
\hat\sigma_4 & \left( \dfrac{\kappa_4 \kappa_5^2}{\kappa_6^3}\right)^{\frac{1}{6}} & \num{5.02e-02}
\\
& & &
\hat\sigma_6 & \left( \dfrac{\kappa_4 \kappa_5^2}{\kappa_6^3}\right)^{\frac{1}{3}}                       & \num{2.52e-03}
\end{tblr}$ \\
\end{talltblr}\\

Note from \cref{tab:sumparams} that the dimensionless parameters $\hat\alpha,\hat\nu_1, \hat\sigma_2,\hat\sigma_3$ and $\hat\sigma_4$ are all correlated with $\hat\sigma_6$. In particular, in the asymptotic limit $\hat\sigma_6\to 0$, all these parameters also tend to zero, although at different speeds.
The glycolytic oscillation persists for any $\hat\sigma_6 > 0$, but the system's geometry degenerates when $\hat\sigma_6 \to 0$; this defines system \eqref{eq:rndxy} as a singular perturbation problem.

\begin{figure}
\centering
\includegraphics[width=0.9\textwidth]{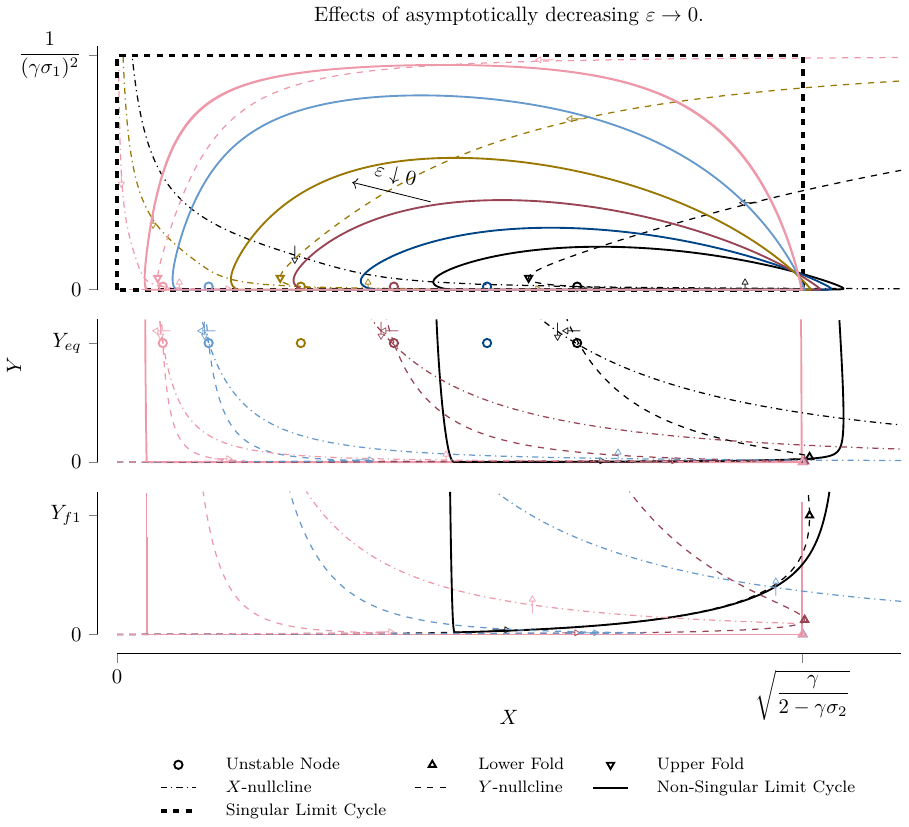}
\caption{Non-dimensionalised System X-Y \eqref{eq:nddynxy}: \emph{Top:} Phase plane showing limit cycles. \emph{Middle:} Zoomed view showing convergence of the equilibrium point $\circ$ at $\mathcal{O}(Y_{eq})$ scale. \emph{Bottom:} Zoomed view showing convergence of the fold point $\triangle$ at $\mathcal{O}(Y_{f1})$ scale. 
}
\label{fig:sysxzntlim}
\end{figure}

\begin{remark}
    To verify that this is a meaningful limit, we demonstrate that it preserves the key geometric features of the oscillator; see also Figure~\ref{fig:sysxzntlim}: 
\begin{enumerate}[label=(\roman*), nosep]
\item the domain is compact and fold points are confined, and
\item there is a unique unstable node.
\end{enumerate}
The first condition is shown to be satisfied in \cref{sec:theoderiv} and the second in \cref{sec:unsnode}.
\end{remark}

\Cref{tab:sumparams} shows that all small parameters are correlated through $\hat\sigma_6$, and we define $\varepsilon := \hat\sigma_6^\frac{1}{4}$. Both $\hat\sigma_3$ and $\hat\sigma_4$ are proportional to $\varepsilon^2$. In fact, $\hat\sigma_4=\varepsilon^2$ based on \cref{tab:sumparams}.
Likewise, $\hat\alpha$ and $\hat\nu_1$ are proportional to $\varepsilon$, giving $\hat\alpha, \hat\nu_1 \sim \mathcal{O}(\varepsilon)$.
The prefactors of the above four small parameters are considered to be of order $\mathcal{O}(1)$.
We have freedom in the exponent $a$ (see \eqref{eq:rels6s2}) of the $\hat\sigma_2$ scaling, and we choose $a=1/2$.
For this scaling $\hat\sigma_2 \sim \mathcal{O}(\sqrt{\hat\sigma_6})$ (i.e., $a=1/2$) to hold, the expression for $\hat\sigma_2$ in \cref{tab:sumparams} implies a constraint on the underlying biophysical parameters $(\kappa_2,\kappa_6)$. Specifically, since  $\hat\sigma_2 =(\kappa_6/\kappa_2)\hat\sigma_6$, the ratio $\kappa_6/\kappa_2$ must scale as $\mathcal{O}(\hat\sigma_6^{-1/2})$ or, equivalently, $\mathcal{O}(\varepsilon^{-2})$.
Altogether, this gives the resulting hierarchy of identified small parameters
\begin{align}
\hat\sigma_6 &= \varepsilon^4 \notag \\
\hat\sigma_2, \hat\sigma_3, \hat\sigma_4 &\sim \mathcal{O}(\varepsilon^2) 
\label{eq:sphier} \noeqref{eq:sphier} \\
\hat\alpha, \hat\nu_1 &\sim \mathcal{O}(\varepsilon) \notag
\end{align}
%
%
With these asymptotic relationships, and dropping the hats from all parameters as well as the subscript in $\hat{\nu}_1$, we transform system \eqref{eq:rndxy} to the singularly perturbed system: 
\begin{align}
\begin{rcases}
\dfrac{dX}{d \tau} &= \varepsilon \left( \alpha - \nu \dfrac{X^2 Y + \varepsilon^4 X^2}{\sigma_1 X^2 Y + \varepsilon^2 \left( \sigma_2 X^2 + \sigma_3 Y + 1 \right)} \right) \\[1em] 
\dfrac{dY}{d \tau} &=  \dfrac{X^2 Y + \varepsilon^4 X^2}{\sigma_1 X^2Y+ \varepsilon^2 \left( \sigma_2 X^2 + \sigma_3 Y + 1 \right) } - \gamma \sqrt{Y} 
\end{rcases} \label{eq:nddynxy} \noeqref{eq:nddynxy}
\end{align}

The singular nature of the oscillation is revealed in \cref{fig:sysxzntlim}, where we simulated system \eqref{eq:nddynxy} with decreasing values of $\varepsilon$. In the limit $\varepsilon\to 0$, the limit cycle approaches the rectangular box (dashed black box), with the $\varepsilon$-dependent frequency tending to zero, indicating that the period grows without bound. The rectangular box also represents the compact domain that emerged from the asymptotic boundedness conditions we imposed. In the middle panel of \cref{fig:sysxzntlim}, we further observe that the unstable node ($\circ$) converges to the point $(0,Y_{eq})$ on the $Y$-axis, where $Y_{eq} \sim \mathcal{O}(1)$. In contrast, the bottom panel of \cref{fig:sysxzntlim} shows that the lower fold point $(X_{f1}, Y_{f1})$ converges to a point on the $X$-axis, as its coordinates approach $(X_{f1,0}, 0)$ where $X_{f1,0} \sim \mathcal{O}(1)$. Both of these observations highlight the singular nature of the oscillation.

\subsection{Eliminating the square root}
In system \eqref{eq:nddynxy}, we notice that $Y=0$ constitutes a singularity due to the $\sqrt{Y}$ term with its infinite slope at $Y=0$. We use a two-step transformation process to deal with this singularity.
In a first step, we perform the coordinate transformation $Y = Z^2$ to obtain a (now) rational vector field

\begin{equation}
\begin{rcases}
\dfrac{dX}{d \tau} &= \varepsilon \left( \alpha - \nu \dfrac{X^2 Z^2 + \varepsilon^4 X^2}{\sigma_1 X^2 Z^2 + \varepsilon^2 ( \sigma_2 X^2 + \sigma_3 Z^2 + 1)} \right) \\[1em] 
\dfrac{dZ}{d \tau} &=  \dfrac{1}{2Z} \left( \dfrac{X^2 Z^2 + \varepsilon^4 X^2}{\sigma_1 X^2 Z^2+ \varepsilon^2 ( \sigma_2 X^2 + \sigma_3 Z^2 + 1) } - \gamma Z \right) 
\end{rcases} \label{eq:nddynxz} \noeqref{eq:nddynxz}
\end{equation}
which exposes the singularity at $Z=0$ in its denominator.

In a second step, we make a state-dependent (nonlinear) time transformation\footnote{This transformation is, by its nature, singular at $Z=0$.} 
\begin{equation}
d\tau = 2Z \left( \sigma_1 X^2 Z^2+ \varepsilon^2 ( \sigma_2 X^2 + \sigma_3 Z^2 + 1) \right) d\bar{\tau} \label{eq:nltimtrans} \noeqref{eq:nltimtrans} \\
\end{equation}
to obtain a polynomial vector field from system \eqref{eq:nddynxz}:

\begin{align}
&\begin{rcases}
X' &= -2 \varepsilon Z \left( (\nu -\alpha \sigma_1) X^2 Z^2 - \varepsilon^2 \alpha (\sigma_2 X^2 + \sigma_3 Z^2 + 1)  + \varepsilon^4 \nu X^2 \right)\\ 
Z' &= (1 - \gamma \sigma_1 Z) X^2 Z^2  - \varepsilon^2 \gamma Z (\sigma_2 X^2 + \sigma_3 Z^2 + 1) + \varepsilon^4 X^2 
\end{rcases} \label{eq:NdSysXZ} \noeqref{eq:NdSysXZ}
\end{align}
where $'$ denotes differentiation with respect to $\bar{\tau}$. Note that the singularity at $Z=0$ of the $dZ/d\tau$-equation in \eqref{eq:nddynxz} has now become an additional nullcline of the $X'$-equation in \eqref{eq:NdSysXZ}.

\begin{figure}
\centering
\includegraphics[width=0.8\textwidth]{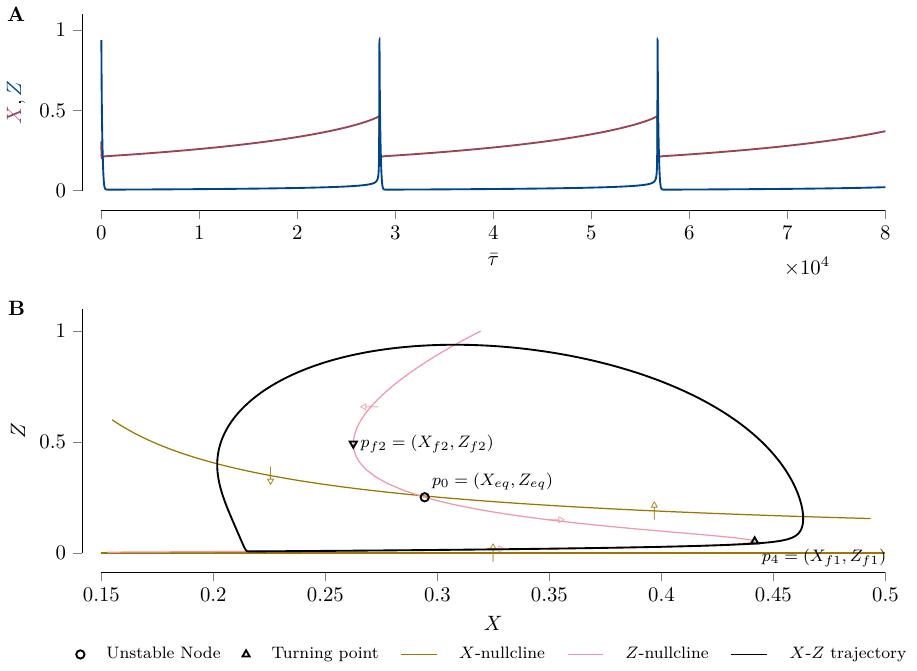}
\caption{Relaxation Oscillator: System X-Z \eqref{eq:NdSysXZ}: \emph{Upper:} time trace, \emph{Lower:} phase plot.}
\label{fig:sysxzscn}
\end{figure}

As expected, this new system \eqref{eq:NdSysXZ} no longer has the two-stroke structure comprising refill and pulse phases as in \cref{fig:origprob}, but rather exhibits relaxation oscillator-type dynamics evolving on multiple timescales (\cref{fig:sysxzscn}), albeit in a degenerate/non-standard form, as will become clear in the following analysis. 
Note, the concomitant stretching of the slow `refill' phase and the shrinkage of the fast `pulse' phase are a direct consequence of the state-dependent time transformation \eqref{eq:nltimtrans}. 

We have thus finished the pre-processing of the Smolen model, arriving at a non-dimensionalised system \eqref{eq:NdSysXZ} with a singular perturbation parameter $\varepsilon\ll 1$ that is amenable to GSPT analysis. 
The key geometric features collapse onto the coordinate axes:\\

\begin{itemize}[left=0pt]
    \item \textit{The unstable node $p_0$:} $(X_{eq}, Z_{eq}) \sim (\mathcal{O}(\varepsilon), \mathcal{O}(1))$. It converges to $(0, Z_{eq,0})$ on the $Z$-axis, where $Z_{eq,0} = \alpha/(\gamma\nu)$ (from \eqref{eq:eqbmxy}).
    \item \textit{The lower fold $p_4$:} $(X_{f1}, Z_{f1}) \sim (\mathcal{O}(1), \mathcal{O}(\varepsilon^2))$. It converges to $(X_{f1,0}, 0)$ on the $X$-axis. Its $\mathcal{O}(1)$ coordinate is $X_{f1,0}^2 = \gamma / (2 - \gamma \sigma_2)$, as derived from \eqref{eq:xf1} and the final hierarchy in \eqref{eq:sphier} (where $\hat\sigma_4 = \varepsilon^2$).
    \item \textit{The upper fold $p_{f2}$:} $(X_{f2}, Z_{f2}) \sim (\mathcal{O}(\varepsilon), \mathcal{O}(1))$. It also converges to the $Z$-axis at $(0, Z_{f2,0})$, where $Z_{f2,0} = \sqrt{Y_{f2,0}} \sim \mathcal{O}(1)$ is the explicit $\mathcal{O}(1)$ root given by \eqref{Yf2-root}.\\
\end{itemize}

\begin{figure}[ht]
\centering
\includegraphics[width=0.8\textwidth]{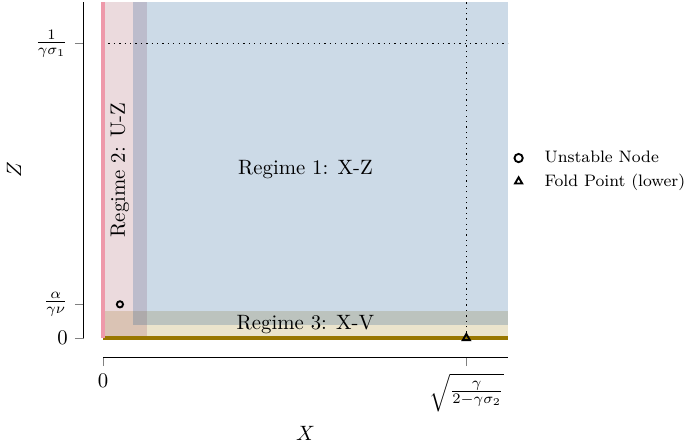}
\caption{The three scaling regimes required for the GSPT analysis, as derived from the asymptotic scalings of the equilibrium and fold points.}
\label{fig:anasysreg}
\end{figure}

This collapse is also seen in the nullclines. In the singular limit, the $X'$-nullcline (\eqref{eq:NdSysXZ}, $X'=0$) degenerates to the axes $Z=0$ and $X=0$. The $Z'$-nullcline ($Z'=0$) degenerates to $Z=0$, $X=0$, and the line $Z=1/(\gamma\sigma_1)$. The extensive overlap on the axes highlights the degeneracy of the singular problem.
\noindent
This multi-scale singular structure, where features exist at $\mathcal{O}(1)$, $\mathcal{O}(\varepsilon)$ and $\mathcal{O}(\varepsilon^2)$ scales, dictates our GSPT analysis strategy. We must analyse the system in three distinct scaling regimes, as illustrated in \cref{fig:anasysreg}:
\begin{itemize}
    \item \textit{Regime 1} ($X, Z \sim \mathcal{O}(1)$): The `outer' view.
    \item \textit{Regime 2} ($X = \varepsilon U, Z \sim \mathcal{O}(1)$): A `zoom' near the $Z$-axis to resolve the unstable node.
    \item \textit{Regime 3} ($X \sim \mathcal{O}(1), Z = \varepsilon^2 V$): A `zoom' near the $X$-axis to resolve the lower fold.\\
\end{itemize}

\section{ {GSPT analysis} in the different scaling regimes}
\label{sec:ana}

{Before detailing the dynamics within the specific scaling regimes, we note that the GSPT analysis throughout this section is carried out in the general formulation \cite{Wechselberger2020}. While some of our derived systems are presented in the standard fast-slow form, we consistently adopt the general framework for two primary reasons. First, this representation highlights the distinct processes evolving on multiple separate $\mathcal{O}(\varepsilon^i)$ timescales. Second, the general formulation is coordinate-independent, allowing us to identify the governing geometric objects, such as the critical manifold, the layer problem, and the reduced flow, without requiring an a priori separation of the phase space variables into fast and slow components.}

\subsection{Regime 1: The X-Z system}\label{sec:regime1}

{We first note that system \eqref{eq:NdSysXZ} is in standard form, with $Z$ as the fast variable and $X$ as the slow one. However, applying the general form of GSPT (see \cref{app:gspt} for further details) gives:}

\begin{equation}
\begin{aligned}
\begin{pmatrix}
X' \\
Z'
\end{pmatrix}
&= 
\begin{pmatrix}
0\\
1
\end{pmatrix}
(1 - \gamma \sigma_1 Z) X^2 Z^2
+ \varepsilon
\begin{pmatrix}
-2  (\nu -\alpha \sigma_1) X^2 Z^3\\
0 
\end{pmatrix}\\
&+ \varepsilon^2
\begin{pmatrix}
0 \\
- \gamma Z (\sigma_2 X^2  + \sigma_3 Z^2 + 1) 
\end{pmatrix}
+ \varepsilon^3
\begin{pmatrix}
2  \alpha Z (\sigma_2 X^2 +  \sigma_3 Z^2 + 1)\\
0
\end{pmatrix}\\
&+ \varepsilon^4
\begin{pmatrix}
0\\
X^2
\end{pmatrix}
+ \varepsilon^5
\begin{pmatrix}
- 2  \nu X^2 Z\\
0
\end{pmatrix} 
:= N_0(X,Z) f_0(X,Z) + \sum_{i=1}^5 \varepsilon^i F_i(X,Z), 
\end{aligned} 
\label{eq:spxz} \noeqref{eq:spxz}
\end{equation}
with $N_0(X,Z)=N_0=(0,1)^\top$ a constant vector and $f_0(X,Z)=(1 - \gamma \sigma_1 Z) X^2 Z^2$ a scalar function. This form highlights the 6 distinct processes evolving on distinct timescales of $\mathcal{O}(\varepsilon^i),\,i=0,1,\ldots 5$.

\subsubsection{The layer problem}

\begin{figure} 
\centering
\includegraphics[width=0.8\textwidth]{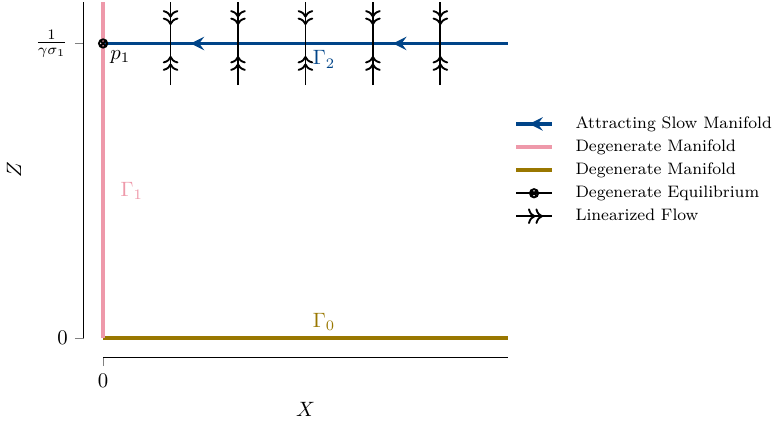}
\caption{Analysis of the X-Z System (Regime 1), showing the three branches of the critical manifold $S_0$ and the reduced flow on the attracting branch $\Gamma_2^N$.}
\label{fig:anasysxz} 
\end{figure}

The layer problem (see \cref{app:gspt}) is obtained by taking the singular limit $\varepsilon \to 0$ in \eqref{eq:spxz}, which gives
\begin{equation}
\begin{pmatrix}
X' \\
Z'
\end{pmatrix}
= 
\begin{pmatrix}
0 \\
1
\end{pmatrix}
(1 - \gamma \sigma_1 Z) X^2 Z^2 =N_0 f_0(X,Z) \,.\label{eq:xzlp} \noeqref{eq:xzlp}
\end{equation}
This is a standard singular perturbation problem: $X$ is a constant (and thus identified as a slow variable), $Z$ is the `fast' variable, and the span of $N_0$ represents the fast fibre of the layer problem. The set of equilibria $S_0$ of the layer problem is defined by
\begin{equation*}
f_0 (X,Z) = X^2 Z^2 (1-\gamma \sigma_1 Z) =0
\end{equation*}
comprising three one-dimensional branches, i.e.,
\begin{equation}
S_0 = \Gamma_0 \cup \Gamma_1 \cup \Gamma_2 = \{Z=0\} \cup \{X=0\} \cup \{Z=\frac{1}{\gamma \sigma_1}\} \label{eq:mnfdsxz} \noeqref{eq:mnfdsxz}  
\end{equation}
that pairwise intersect at $\Gamma_0 \cap \Gamma_1=(0,0)$ and $p_1 := \Gamma_1 \cap \Gamma_2=(0,1/(\gamma\sigma_1))$; see Figure \ref{fig:anasysxz}. Each branch is viewed as a one-dimensional critical manifold of the layer problem. To determine the stability properties of each branch $\Gamma_i,\,i=0,1,2$, we calculate the non-trivial eigenvalue:
\begin{equation}
\lambda = Df_0 N_0 = D_Z f_0 = X^2 Z(2  - 3\gamma \sigma_1 Z)\,. \label{eq:xzdegm}
\end{equation}
We notice that this eigenvalue vanishes when evaluated along the branches $\Gamma_0$ and $\Gamma_1$, and, hence, both branches are degenerate. On the other hand, $\lambda$ evaluated along $\Gamma_2$ gives
\begin{equation*}
\lambda(X,\frac{1}{\gamma\sigma_1}) = -\frac{X^2}{\gamma \sigma_1} \leq 0.
\end{equation*}
Hence, we have $\Gamma_2=\Gamma_2^{N-}\cup \{p_1 \} \cup \Gamma_2^{N+}$ with normally hyperbolic branches $\Gamma_2^{N\pm}$, both attracting, and loss of normal hyperbolicity at $p_1$.

\subsubsection{The reduced problem along $\Gamma_2$}\label{sec:regime1reduced}

The corresponding leading order reduced problem (see \cref{app:gspt}) along $\Gamma_2$ is given by 
\begin{equation}
\dot{X} = -2  (\nu - \alpha \sigma_1) X^2 Z^3|_{\Gamma_2}=
-2  \frac{(\nu - \alpha \sigma_1)}{\gamma^3 \sigma_1^3} X^2\,.
\label{eq:xzrp} \noeqref{eq:xzrp}
\end{equation}
Since $\nu > \alpha \sigma_1$, see \eqref{cond:param2}, the flow on the slow manifold $\Gamma_2$ decreases monotonically towards $X = 0$ for $X>0$.

\subsection{Regime 2: The U-Z system}\label{sec:regime2}

To understand the slow dynamics near the degenerate manifold $\Gamma_1$ of the X-Z system \eqref{eq:spxz}, we zoom in on Regime 2 through the following coordinate and time transformations:
\begin{equation}
X := \varepsilon U, \quad d\bar{\tau} = \varepsilon^{-2} d\tau_2\,.
\label{eq:uzdt2} \noeqref{eq:uzdt2}
\end{equation}
With $' = d/d\tau_2$ and $\tilde{\varepsilon} := \varepsilon^2$, we obtain
\begin{equation}
\begin{rcases}
U' &= -2 Z \left[ (\nu - \alpha \sigma_1) U^2 Z^2 - \alpha(\sigma_3 Z^2 + 1) - \tilde{\varepsilon} \alpha \sigma_2 U^2 + \tilde{\varepsilon}^2 \nu U^2  \right] \\
Z' &= (1- \gamma \sigma_1 Z) U^2 Z^2 - \gamma Z (\sigma_3 Z^2 + 1) - \tilde{\varepsilon} \gamma \sigma_2 U^2 Z + \tilde{\varepsilon}^2 U^2
\end{rcases} \label{eq:nduz} \noeqref{eq:nduz}
\end{equation}
which is a singular perturbation problem but in general, \emph{non-standard form}\footnote{In non-standard form, there is no explicit identification of slow and fast variables.} 
\begin{equation}
\begin{aligned}
\begin{pmatrix}
U' \\
Z'
\end{pmatrix}
& = 
\begin{pmatrix}
-2 \left[ (\nu - \alpha \sigma_1) U^2 Z^2 - \alpha(\sigma_3 Z^2 + 1) \right] \\
(1- \gamma \sigma_1 Z) U^2 Z - \gamma (\sigma_3 Z^2 + 1)
\end{pmatrix}
Z
+ \tilde{\varepsilon}
\begin{pmatrix}
2 \alpha \sigma_2 U^2 Z\\
- \gamma \sigma_2 U^2 Z
\end{pmatrix}
+ \tilde{\varepsilon}^2
\begin{pmatrix}
- 2 \nu U^2 Z  \\
U^2
\end{pmatrix} 
\\
&:= N_0(U,Z) f_0(U,Z) + \sum_{i=1}^2 \tilde{\varepsilon}^i F_i(U,Z),
\end{aligned} \label{eq:spuz} \noeqref{eq:spuz}
\end{equation}
with state-dependent vector $N_0(U,Z)$ and scalar function $f_0(U,Z)=Z$. This form highlights the 3 distinct processes evolving on the timescales $\mathcal{O}(\tilde\varepsilon^i),\, i =0,1,2$. 

\subsubsection{The layer problem}
The layer problem is obtained by taking the singular limit $\tilde{\varepsilon} \to 0$ in \eqref{eq:spuz}, which gives
\begin{equation}
\begin{pmatrix}
U' \\
Z'
\end{pmatrix}
= 
\begin{pmatrix}
-2 \left[ (\nu - \alpha \sigma_1) U^2 Z^2 - \alpha(\sigma_3 Z^2 + 1) \right] \\
(1- \gamma \sigma_1 Z) U^2 Z - \gamma (\sigma_3 Z^2 + 1)
\end{pmatrix}
Z = N_0(U,Z) f_0(U,Z).
\label{eq:uzlp} \noeqref{eq:uzlp} 
\end{equation}
The set of equilibria $S_0$ of \eqref{eq:uzlp} comprises the one-dimensional critical manifold $\Gamma_5 = \{Z=0 \}$ and the single isolated point $p_0$ at which the vector $N_0(U,Z)$ vanishes. For the critical manifold $\Gamma_5$, the Jacobian is given by
\begin{equation*}
J(\Gamma_5) = N_0(\Gamma_5) Df_0(\Gamma_5) =
\begin{pmatrix}
0 & 2 \alpha \\
0 &  - \gamma 
\end{pmatrix}
\end{equation*}
with the non-trivial eigenvalue $\lambda(\Gamma_5) = -\gamma < 0$, indicating that the critical manifold $\Gamma_5$ is uniformly attracting along the eigendirection $[2 \alpha,  -\gamma]^\top$.
To analyse the stability of the equilibrium point $p_0$, we evaluate the Jacobian:
\begin{equation*}
\tilde{J}(U,Z) = D{N_0}(U,Z)f_0(U,Z)
= 
\begin{pmatrix}
-4 (\nu - \alpha \sigma_1) U Z^3 & -4 ((\nu - \alpha \sigma_1) U^2  - \alpha \sigma_3 )Z^2\\
2 (1 - \gamma \sigma_1 Z) U Z^2 & ((1 - 2 \gamma \sigma_1 Z) U^2  - 2 \gamma \sigma_3 Z)Z
\end{pmatrix}
\label{eq:sysuzjacp0}
\end{equation*}
with trace $\tr(\tilde{J}(U,Z)) = -4 (\nu - \alpha \sigma_1) U Z^3 + \left[ (1 - 2 \gamma \sigma_1 Z) U^2 - 2 \gamma \sigma_3 Z \right] Z$ which expands to:
$$\tr(\tilde{J}(U,Z)) = -4\nu U Z^3 + 4\alpha\sigma_1 U Z^3 + U^2 Z - 2\gamma\sigma_1 U^2 Z^2 - 2\gamma\sigma_3 Z^2\,.$$
Since $\det(\tilde{J}(p_0)) > 0$, $\tr(\tilde{J}(p_0)) > 0$ and $\tr(\tilde{J}(p_0))^2 > 4 \det(\tilde{J}(p_0))$, we conclude that $p_0$ is an unstable node (see also \cref{app:sec:linana}).

\begin{figure}
\centering
\includegraphics[width=0.9\textwidth]{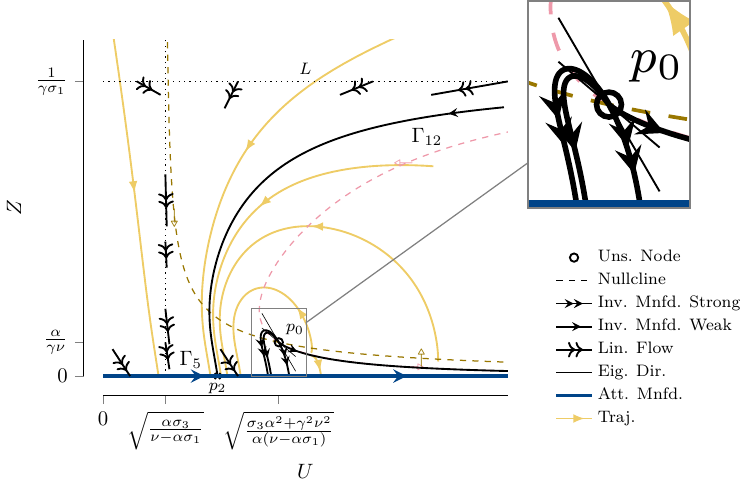}
\caption{Analysis of the U-Z System (Regime 2). The layer problem \eqref{eq:uzlp} shows the phase space partitioned by the $U$-nullcline (olive), $Z$-nullcline (magenta), and a separatrix $\Gamma_{12}$. The flow is dominated by the unstable node $p_0$ and the attracting manifold $\Gamma_5 = \{Z=0\}$.}
\label{fig:anasysuz}
\end{figure} 

The layer flow, shown in \cref{fig:anasysuz}, is determined by the interplay between the unstable node $p_0$ and the attracting critical manifold $\Gamma_5$. The node $p_0$ repels all nearby trajectories, which are then drawn towards the uniformly attracting manifold $\Gamma_5 = \{Z=0\}$. Overall, all trajectories (except the node $p_0$ itself) eventually approach this attracting manifold, which provides the initial condition for the slow flow.

\begin{remark}
   While the local layer flow simply approaches $\Gamma_5$, the nullclines and the separatrix $\Gamma_{12}$ shown in \cref{fig:anasysuz} are important for understanding the global dynamics. As a later analysis will show, $\Gamma_{12}$ (\cref{fig:3dchartk1}) is the true extension of the slow flow from the manifold $\Gamma_2$ and, hence, also plays an important role in partitioning trajectories arriving from infinity. 
\end{remark}

\subsubsection{The reduced problem}

To obtain the slow flow, we introduce the slow timescale $\tilde\tau_s$ through $d\tilde\tau_s = \tilde\varepsilon d\tau_2$, which transforms the fast problem \eqref{eq:spuz} into the slow problem
\begin{equation}
\begin{aligned}
\tilde\varepsilon
\begin{pmatrix}
\dot{U} \\
\dot{Z}
\end{pmatrix}
&= 
\underbrace{
\begin{pmatrix}
-2 \left[ (\nu - \alpha \sigma_1) U^2 Z^2 - \alpha(\sigma_3 Z^2 + 1) \right] \\
(1- \gamma \sigma_1 Z) U^2 Z - \gamma (\sigma_3 Z^2 + 1)
\end{pmatrix}
Z
}_{F_0(U,Z)=N_0(U,Z)f_0(U,Z)}
+ \tilde{\varepsilon}
\underbrace{
\begin{pmatrix}
2 \alpha \sigma_2 U^2 Z\\
- \gamma \sigma_2 U^2 Z
\end{pmatrix}
}_{F_1(U,Z)}
+ \tilde{\varepsilon}^2
\underbrace{
\begin{pmatrix}
- 2 \nu U^2 Z  \\
U^2
\end{pmatrix}
}_{F_2(U,Z)}
=F(U,Z,\tilde\varepsilon)
\end{aligned}
\end{equation}
where $\dot{} = d/d\tilde\tau_s$. In the singular limit $\tilde\varepsilon \to 0$, this system is only well defined along the critical manifold $\Gamma_5 = \{Z=0\}$.

{Along $\Gamma_5$, the reduced flow is obtained by projecting the $\mathcal{O}(\tilde\varepsilon)$ vector field $F_1$ onto the tangent space of $\Gamma_5$, see \cref{fig:app:gsptpm}. Evaluating $F_1$ on $\Gamma_5$ (where $Z=0$) gives $F_1(U,0)=[0,0]^\top$, so this projection vanishes and the drift along $\Gamma_5$ is not determined at $\mathcal{O}(\tilde\varepsilon)$. 
The dynamics along $\Gamma_5$ are therefore governed by higher-order terms, and obtaining them requires an accurate description of the invariant manifold itself. We thus employ the parametrisation method \cite{Lizarraga2021}, which constructs the manifold and the flow on it simultaneously as formal power series in $\tilde\varepsilon$, order by order (see \cref{app:pm} for the detailed mathematical formulation).

Applying this method, we find that the $\mathcal{O}(\tilde\varepsilon)$ corrections to both the manifold embedding and the slow vector field are trivial. Advancing to $\mathcal{O}(\tilde\varepsilon^2)$, and evaluating the non-trivial contribution $F_2$ on the unperturbed manifold, we obtain the leading-order non-trivial reduced vector field (see \cref{app:pm} for the derivation):
\begin{equation}
r_2(\xi) = \frac{2 \alpha \xi^2}{\gamma}\ge 0\,.
\label{eq:uv:sofl} \noeqref{eq:uv:sofl}
\end{equation}
where the coordinate $\xi=U$ parametrises $\Gamma_5$. Thus, the flow on the slow manifold $\Gamma_5$ is non-trivial at $\mathcal{O}(\tilde{\varepsilon}^2)$ and increases monotonically, moving away from the $Z$-axis.}

\subsubsection{Summary of Regime 2}

The analysis in Regime 2 resolves the dynamics near the degenerate manifold $\Gamma_1$. We find that the unstable node $p_0$ is preserved at $\mathcal{O}(1)$ coordinates in this scaling. The layer flow pushes trajectories away from this node and onto the attracting critical manifold $\Gamma_5 = \{Z=0\}$. The reduced problem on $\Gamma_5$ is non-trivial at $\mathcal{O}(\tilde{\varepsilon}^2)$ and consists of a slow flow $r_2$ moving monotonically to the right (increasing $U$). This flow eventually carries the trajectory out of Regime 2, necessitating a transition to the next regime.

\subsection{Regime 3: The X-V system}\label{sec:regime3}

To understand the dynamics near the degenerate manifold $\Gamma_0$ of the X-Z system \eqref{eq:spxz}, we zoom in on Regime 3. We use the coordinate transformation $Z := \varepsilon^2 V$ and the fast time $\tau_3$ defined by the differential relation:
\begin{equation}
d\bar{\tau} = \varepsilon^{-2} d\tau_3\,.
\label{eq:sysxvtimtrans} \noeqref{eq:sysxvtimtrans}
\end{equation}
With $' = d/d\tau_3$, the original system \eqref{eq:spxz} transforms into:
\begin{equation}
\begin{rcases}
X' &=  2 \varepsilon^3 V \left[ \alpha (\sigma_2 X^2 + 1) - \varepsilon^2  X^2  ((\nu -\alpha \sigma_1)  V^2 + \nu)  + \varepsilon^4 \alpha \sigma_3 V^2 \right] \\
V' &= (1 - \varepsilon^2 \gamma \sigma_1 V) X^2 V^2 - \gamma V (\sigma_2 X^2 + 1) + X^2 - \varepsilon^4 \gamma \sigma_3 V^3
\end{rcases} \label{eq:nduv} \noeqref{eq:nduv}
\end{equation}
{Consistent with our overarching approach, we present system \eqref{eq:nduv} in the general formulation:}
\begin{equation}
\begin{aligned}
\begin{pmatrix}
X' \\
V'
\end{pmatrix}
&= 
\underbrace{
\begin{pmatrix}
0\\
1
\end{pmatrix}
(X^2 V^2 + X^2 - \gamma \sigma_2 X^2 V - \gamma V)
}_{N_0 f_0(X,V)}
+ \varepsilon^2
\underbrace{
\begin{pmatrix}
0\\
-\gamma \sigma_1 X^2 V^3
\end{pmatrix}
}_{F_2(X,V)}
+ \varepsilon^3
\underbrace{
\begin{pmatrix}
2 \alpha (\sigma_2 X^2 + 1) V\\
0 
\end{pmatrix}
}_{F_3(X,V)}\\ 
&+ \varepsilon^4
\begin{pmatrix}
0\\
- \gamma \sigma_3 V^3
\end{pmatrix}
+ \varepsilon^5
\begin{pmatrix}
-2 X^2 V \left( (\nu - \alpha \sigma_1) V^2 + \nu \right) \\ 
0
\end{pmatrix}
+ \varepsilon^7
\begin{pmatrix}
2 \alpha \sigma_3 V^3 \\ 
0
\end{pmatrix}\\
&:= N_0 f_0(X,V) + \sum_{i} \varepsilon^i F_i(X,V), \quad i \in \{2,3,4,5,7\}.
\end{aligned} \label{eq:spuv} \noeqref{eq:spuv}
\end{equation}
to highlight the 6 distinct processes evolving on the timescales $\mathcal{O}(\varepsilon^i),\, i =0,2,3,4,5,7$.

\subsubsection{The layer problem}
\begin{figure}[t]
\centering
\includegraphics[width=0.9\textwidth]{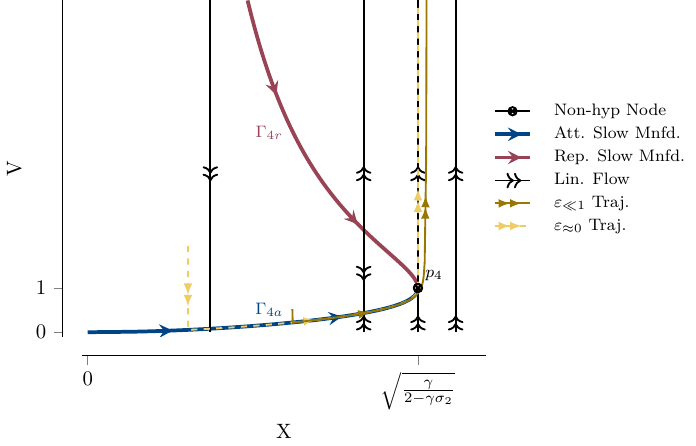}
\caption{Analysis of the X-V System (Regime 3). The phase space dynamics is organised by the critical manifolds $\Gamma_{4a}$ and $\Gamma_{4r}$, attracting and repelling respectively, connected at the non-hyperbolic point $p_4$. The fast trajectories of the layer problem \eqref{eq:uvlp} are attracted towards $\Gamma_{4a}$, onto which they move slowly, according to the reduced problem \eqref{eq:spuv-slow}, until they are discharged at $p_4$, back to the fast layer phase. These are demonstrated by typical trajectories: non-singular (brown) and (nearly) singular (yellow).}  
\label{fig:anasysuv}
\end{figure}

The layer problem is obtained by taking the limit $\varepsilon \to 0$:
\begin{equation}
\begin{pmatrix}
X' \\
V'
\end{pmatrix}
= 
\begin{pmatrix}
0 \\
1
\end{pmatrix}
(X^2 V^2 + X^2 - \gamma V -\gamma \sigma_2 X^2 V)
= N_0 f_0(X,V).  \label{eq:uvlp} \noeqref{eq:uvlp} 
\end{equation}
The set of equilibria $S_0$, defined by $f_0(X,V)=0$, forms the critical manifold. By solving for $X^2$, we find $X^2 = \gamma V / (V^2 - \gamma \sigma_2 V + 1)$. This defines the one-dimensional folded critical manifold $\Gamma_4$:

\begin{equation}\label{Gamm4:para}
\Gamma_4 = \left\{ \left. \left(\sqrt{\frac{\gamma V}{V^2 - \gamma \sigma_2 V + 1}},  V  \right) \right\rvert V \geq 0  \right\}
\end{equation}
The stability of the critical manifold $\Gamma_4$ is determined by the non-trivial eigenvalue $\lambda(X,V) = D_V f_0(X,V) = X^2 (2V - \gamma \sigma_2) - \gamma$, which, evaluated along $\Gamma_4$, gives:
\begin{equation}
\lambda(\Gamma_4) = \frac{\gamma (V^2 - 1)}{V^2 - \gamma \sigma_2 V + 1} 
\begin{cases}
< 0, & V \in \Gamma_{4a} := \{\Gamma_4: 0 \le V < 1 \} \\ 
= 0, & V = 1 \\
> 0, & V \in \Gamma_{4r} := \{\Gamma_4: V > 1 \}
\end{cases}
\end{equation}
Thus, the critical manifold $\Gamma_4 = \Gamma_{4a} \cup \{p_4\} \cup \Gamma_{4r}$ consists of an attracting branch $\Gamma_{4a}$ and a repelling branch $\Gamma_{4r}$ connected at a fold point $p_4 = (\sqrt{\gamma/(2-\gamma\sigma_2)}, 1)$ where normal hyperbolicity is lost; see \cref{fig:anasysuv}.

\subsubsection{The reduced problem}
{
To obtain the slow flow, we introduce the slow timescale $\tau_s$ defined by $d\tau_s = \varepsilon d\tau_3$. This transforms the fast problem \eqref{eq:spuv} into the slow problem
\begin{equation}
\varepsilon
\begin{pmatrix}
\dot{X} \\
\dot{V}
\end{pmatrix}
= N_0 f_0(X,V) + \sum_{i} \varepsilon^i F_i(X,V), \quad i \in \{2,3,4,5,7\} \label{eq:spuv-slow} \noeqref{eq:spuv-slow}
\end{equation}
where $\dot{} = d/d\tau_s$. The singular limit $\varepsilon \to 0$ is only well defined if the phase space is restricted to the critical manifold $\Gamma_4$, i.e., to $f_0(X,V) = 0$.

Along $\Gamma_4$, the leading-order perturbation terms fail to determine the drift. The $\mathcal{O}(\varepsilon)$ vector field $F_1$ is zero. Furthermore, the first non-trivial field $F_2$ is vertical, aligning with the fast fibre $N_0 = [0, 1]^\top$. Since the slow (tangent) space is complementary to this fast fibre, $F_2$ has no component in the slow space, and its projection onto $\Gamma_4$ vanishes.

As in Regime 2, obtaining the non-trivial slow flow requires higher-order terms, i.e., the corrections to the flow and the manifold must be computed simultaneously using the parametrisation method \cite{Lizarraga2021} (see \cref{app:pm} for the detailed derivation). Using $\xi:=V$ to parametrise $\Gamma_4$ as a graph, this method reveals that the $\mathcal{O}(\varepsilon)$ corrections are trivial. At $\mathcal{O}(\varepsilon^2)$, while the reduced vector field remains zero ($r_2 = 0$), $F_2$ induces a non-trivial second-order correction to the manifold embedding:
\begin{equation}
\varphi_2(\xi) = 
\begin{pmatrix}
\dfrac{\sigma_1 \xi^4}{2} \sqrt{\dfrac{\gamma^3}{\xi (\xi^2 - \gamma \sigma_2 \xi +1)^3}}\\
0
\end{pmatrix}
\label{eq:xvpmes} \noeqref{eq:xvpmes}
\end{equation}
Advancing to $\mathcal{O}(\varepsilon^3)$, the combination of $F_3$ (which possesses a non-trivial slow component) and the embedding correction $\varphi_2$ yields the first non-trivial reduced vector field $r_3$:
\begin{equation}
r_3(\xi) = \dfrac{4 \alpha (1+\xi^2)}{1-\xi^2} \sqrt{\dfrac{\xi^3 (\xi^2 - \gamma \sigma_2 \xi +1) }{\gamma}}
\quad \begin{cases}
> 0,  & 0 < \xi < 1 \text{ (on } \Gamma_{4a} \text{)} \\
< 0,  & \xi > 1 \text{ (on } \Gamma_{4r} \text{)}
\end{cases}
\label{eq:xvsf} \noeqref{eq:xvsf}
\end{equation}
Thus, the leading-order slow flow, given by $\dot\xi = \varepsilon^3 r_3(\xi)$, is directed towards the fold point $p_4$ (where $\xi=1$) from both the attracting and repelling branches; see \cref{fig:anasysuv}.

The reduced vector field $r_3(\xi)$ blows up at $\xi = 1$. This is a consequence of the loss of normal hyperbolicity at $p_4$, which is a \emph{regular fold (jump) point} (see \cref{df:app:rjp}). It corresponds to a saddle-node bifurcation in the layer problem, and the blow-up of the reduced flow at this point is its characteristic feature.
}

\subsection{Constructing the singular limit cycle}
\begin{figure}
\centering
\includegraphics[width=0.8\textwidth]{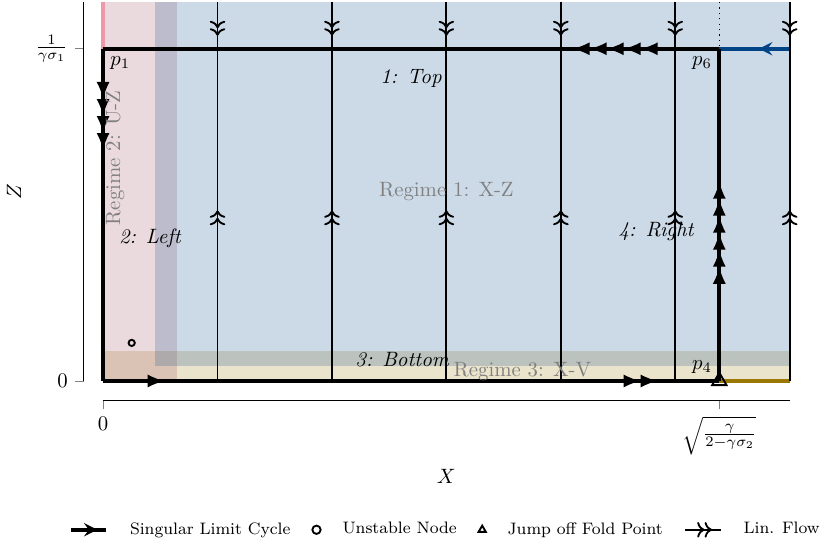}
\caption{
The Singular Limit Cycle. A composite schematic showing the singular manifolds. The cycle consists of four distinct segments: (1) a slow flow on the \emph{Top} hyperbolic manifold from $p_6$ to $p_1$, (2) a transition from $p_1$ to the origin along the degenerate \emph{Left} manifold (governed by Regime 2), (3) a composite ultra-slow/slow flow along the degenerate \emph{Bottom} manifold (governed by Regime 3), and (4) a fast jump from $p_4$ to $p_6$ (through the layer dynamics). Note that the `flows' described on the degenerate \emph {2: Left} and \emph{3: Bottom} manifolds represent the leading-order dynamics revealed by the rescaling analyses in Regimes 2 and 3; see \cref{fig:anasysuz,fig:anasysuv}.}
\label{fig:sysslc}
\end{figure}

Combining the analyses of the different regimes, we can now reconstruct the composite geometry of the singular limit cycle. This reconstruction, shown in \cref{fig:sysslc} (thick black rectangle), is a schematic superposition: it projects the leading-order dynamics revealed by each scaling regime back onto the singular ($\varepsilon=0$) skeleton. In this limit, the degenerate manifolds $\Gamma_1$ and $\Gamma_0$ (containing $\Gamma_5$) are sets of non-hyperbolic equilibria; the `flow' we describe on them is the $\varepsilon$-dependent dynamics that emerges upon rescaling. The singular limit cycle, shown in the X-Z coordinates in \cref{fig:sysslc}, consists of four distinct segments. We trace the cycle starting from $p_6$:
\begin{enumerate}
    \item \textit{Top segment} ($p_6 \to p_1$): This segment is the attracting critical manifold $\Gamma_2$ from \textit{Regime 1}. The flow is governed by the $\mathcal{O}(\varepsilon)$ reduced problem ($\dot{X} < 0$, \eqref{eq:xzrp}), moving trajectories slowly towards the $Z$-axis ($p_1$).

    \item \textit{Left segment} ($p_1 \to \text{origin}$): This segment is governed by \textit{Regime 2}. The trajectory arrives at $p_1$ (on the $Z$-axis, $\Gamma_1$) and follows the \textit{Regime 2 layer problem} ($d/d\tau_2$ time) down the $Z$-axis. This flow ($Z' < 0$, \eqref{eq:uzlp}) guides the trajectory to the origin.

    \item \textit{Bottom segment} ($\text{origin} \to p_4$): This is a composite segment that begins at the origin:
    \begin{itemize}
        \item The trajectory first moves along the normally hyperbolic manifold $\Gamma_5 = \{Z=0\}$ (the $X$-axis), governed by the \textit{Regime 2} slowest \textit{degenerate \footnote{Here, degenerate refers to the dynamics being slower than leading order} reduced problem.}. This $\mathcal{O}(\varepsilon^8)$ flow ($r_2 > 0$, \eqref{eq:uv:sofl}) carries the trajectory away from the origin.
        \item This path is continuous with the attracting branch $\Gamma_{4a}$ of the folded manifold (from \textit{Regime 3}), which also starts at the origin. As the trajectory moves into the $X \sim \mathcal{O}(1)$ region, its dynamics are characterised by the \textit{Regime 3} \textit{degenerate reduced problem}, with a ramped up $\mathcal{O}(\varepsilon^6)$ flow ($r_3 > 0$, \eqref{eq:xvsf}), which accelerates significantly as it approaches the regular fold point $p_4$.
    \end{itemize}

\item \textit{Right segment} ($p_4 \to p_6$): This is a fast jump. Having reached the fold point $p_4$ (defined by the \textit{Regime 3} analysis), the trajectory is discharged from the critical manifold. This jump is initially governed by the \textit{Regime 3 layer problem} (see \eqref{eq:uvlp}) with speed $\mathcal{O}(\varepsilon^2)$, consistent with the Regime 3 layer timescale $d\tau_3$. As the trajectory moves away from the degenerate region ($Z \approx 0$) and $V$ becomes large, the flow accelerates further to the $\mathcal{O}(1)$ dynamics of the \textit{Regime 1 layer problem} to complete the connection to the attracting manifold $\Gamma_2$ at point $p_6$.
\end{enumerate}
\cref{tab:fsssum} summarises the timescales of the flow on each segment, normalised to the original timescale $\bar{\tau}$.\\

\begin{talltblr}[
  caption = {Summary of timescales for each segment of the singular limit cycle, in order of flow.}, %
  label   = {tab:fsssum},                                         %
]{
  colspec = {c} 
}
$\begin{tblr}{
  colspec = {l|l|l|l},
  columns = {mode=math} %
}
\text{Segment} & \text{Manifold} & \text{Analysis} & \text{Flow Speed } (d/d\bar{\tau})  \\
\hline
\text{1. Top} & \Gamma_2 & \text{Regime 1 (Reduced)} & \mathcal{O}(\varepsilon) \text{ (slow)} \\
\text{2. Left} & \Gamma_1 & \text{Regime 2 (Layer)} & \mathcal{O}(\varepsilon^2) \text{ (slower)} \\
\text{3a. Bottom} & \Gamma_5 & \text{Regime 2 (Reduced)} & \mathcal{O}(\varepsilon^8) \text{ (slowest)} \\
\text{3b. Bottom} & \Gamma_{4a} & \text{Regime 3 (Reduced)} & \mathcal{O}(\varepsilon^6) \text{ (infra-slow)} \\
\text{4a. Right} & \text{Layer Flow} & \text{Regime 3 (Layer)} & \mathcal{O}(\varepsilon^2) \text{ (slower)} \\
\text{4b. Right} & \text{Layer Flow} & \text{Regime 1 (Layer)} & \mathcal{O}(1) \text{ (fast)} 
\end{tblr}$ \\
\end{talltblr}

\subsubsection*{Explanation of timescales}
The hierarchy of speeds in \cref{tab:fsssum} is derived by normalising all dynamics to the reference timescale $d\bar{\tau}$ of System \eqref{eq:NdSysXZ}.
\begin{itemize}[left=0pt]
    \item \textit{$\mathcal{O}(1)$ (fast):} This is the reference fast timescale of the \textit{Regime 1 layer problem} \eqref{eq:xzlp}. As required by Fenichel theory, the fast jump from $p_4$ follows these $\mathcal{O}(1)$ dynamics to connect to the attracting manifold $\Gamma_2$ at $p_6$.

    \item \textit{$\mathcal{O}(\varepsilon)$ (slow):} This speed arises from the first non-trivial reduced problem in Regime 1. The flow is driven by the $\mathcal{O}(\varepsilon)$ term in System \eqref{eq:spxz} and governs the slow drift along the normally hyperbolic manifold $\Gamma_2$.

    \item \textit{$\mathcal{O}(\varepsilon^2)$ (slower):} This speed appears in two distinct segments.
    \begin{enumerate}
        \item[(a)] It governs the \textit{layer dynamics} of Regime 2 (left segment, $\Gamma_1$). The local fast time is $d\tau_2 = \varepsilon^2 d\bar{\tau}$ (from \eqref{eq:uzdt2}), so the $\mathcal{O}(1)$ local flow becomes $\mathcal{O}(\varepsilon^2)$ when normalised.
        \item[(b)] It also informs the \textit{layer dynamics} of Regime 3 (right segment), with local fast timescale $d\tau_3 = \varepsilon^2 d\bar\tau$ (from \eqref{eq:sysxvtimtrans}), resulting in a normalised speed of $\mathcal{O}(\varepsilon^2)$.
    \end{enumerate}

    \item \textit{$\mathcal{O}(\varepsilon^6)$ (infra-slow):} This speed is captured in the degenerate reduced problem of Regime 3 (bottom segment, $\Gamma_{4a}$). The local slow time is $d\tau_s = \varepsilon d\tau_3 = \varepsilon^3 d\bar\tau$. In its local timescale, the reduced flow $r_3$ \eqref{eq:xvsf} is of order $\mathcal{O}(\varepsilon^3)$, which gives a normalised speed of $\mathcal{O}(\varepsilon^3) \times \varepsilon^3 = \mathcal{O}(\varepsilon^6)$. 

    \item \textit{$\mathcal{O}(\varepsilon^8)$ (slowest):} This slowest speed is found in the degenerate reduced problem of Regime 2 (bottom segment, $\Gamma_5$). The local slow time is $d\tilde\tau_s = \tilde\varepsilon d\tau_2 = \varepsilon^4 d\bar{\tau}$. The reduced flow $r_2$ is of order $\mathcal{O}(\tilde\varepsilon^2) = \mathcal{O}(\varepsilon^4)$ in this local timescale (from \eqref{eq:uv:sofl}). Normalising this gives a speed of $\mathcal{O}(\varepsilon^4) \times \varepsilon^4 = \mathcal{O}(\varepsilon^8)$, explaining the `crawl' away from the origin.
\end{itemize}
This multi-scale patchwork, comprising various timescale changes across different regimes, is a direct consequence of the singular geometry at the origin, which is rigorously unified using the blow-up method \cite{Dumortier1996, Krupa2001, KrupaW2010, KrupaSz2001, Kosiuk2011, JelbartW2022}.

From \cref{tab:fsssum}, we can deduce the dynamics of the cycle`s speed. The cycle is dominated by its slow ($\mathcal{O}(\varepsilon), \mathcal{O}(\varepsilon^2)$) and crawling ($\mathcal{O}(\varepsilon^8), \mathcal{O}(\varepsilon^6)$) segments. The cycle is `closed' by a fast jump in the layer flow ($\mathcal{O}(\varepsilon^2), \mathcal{O}(1)$) back to the top. The flow `crawls' away from the origin (speed $\mathcal{O}(\varepsilon^8)$), ramps up to $\mathcal{O}(\varepsilon^6)$, followed by an acceleration along the bottom as it approaches $p_4$. It jumps with an initial speed of $\mathcal{O}(\varepsilon^2)$ with significant acceleration to $\mathcal{O}(1)$, relaxes along the top ($\mathcal{O}(\varepsilon)$), and transitions down the $Z$-axis ($\mathcal{O}(\varepsilon^2)$).

We also point out that our non-dimensionalised model \eqref{eq:NdSysXZ} is more complicated than the one inherited and analysed in \cite{Kosiuk2011} which described a single slow vector field for each regime that invariably yielded a standard fast-slow problem. The mechanistic complexity inherent in the original Smolen-based biophysical model \cite{Marinelli2018, Bertram2017, Bertram2023} percolates to the dynamics of the non-dimensionalised system, which we have uncovered here. 

The above three Regimes 1, 2 and 3 have helped decipher the dynamics at different timescales. The above reconstruction of the singular limit is but an artificial superposition of the regions they define, as shown in \cref{fig:anasysreg,fig:sysslc}, since we have not established how the different patches connect, including smooth transitions between timescales. To address this gap, we closely follow the blow-up analysis in \cite{Kosiuk2011} to show how the flows transition between timescales.


\section{Blow-up analysis}\label{sec:blowup}

While the scaling analysis revealed the local dynamics within specific regimes, it did not address the global connectivity between them. To prove the existence of a singular limit cycle and establish the geometric `gluing' of these disparate timescales, we employ the blow-up method \cite{Krupa2001, KrupaSz2001, Kosiuk2011, Dumortier1996}. This technique desingularises the vector field along the degenerate coordinate axes, unfolding non-hyperbolic critical manifolds onto higher-dimensional spaces where partial hyperbolicity can be recovered. Following the approach in \cite{Kosiuk2011}, we perform a sequential blow-up to rigorously verify the heteroclinic connections between the flows of Regimes 1, 2 and 3.

\subsection{First cylindrical blow-up}

\begin{figure}[t]
\centering
\includegraphics[width=0.85\textwidth]{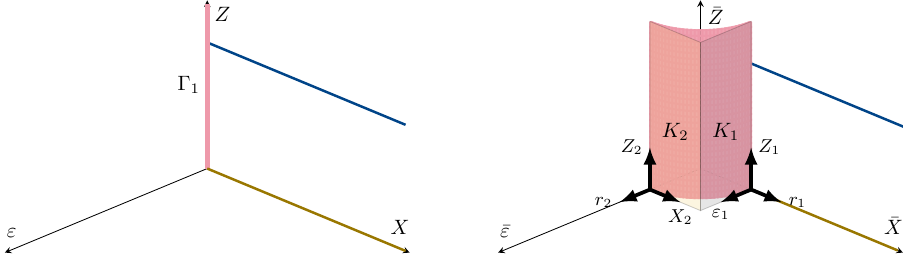}
\caption{First cylindrical blow-up.}
\label{fig:busc1_12}
\end{figure}

We start by performing a first blow-up of System X-Z \eqref{eq:NdSysXZ}, where the two-dimensional vector field is extended to a three-dimensional system by adding the trivial equation $\varepsilon' = 0$:
\begin{align}
\begin{rcases}
X' &= - 2 \varepsilon Z  \left[ (\nu - \alpha \sigma_1) X^2 Z^2 - \varepsilon^2 \alpha (\sigma_2 X^2 + \sigma_3 Z^2 + 1 ) + \varepsilon^4 \nu X^2 \right]\\
Z' &= (1 - \gamma \sigma_1 Z) X^2 Z^2 - \varepsilon^2 \gamma Z (\sigma_2 X^2 + \sigma_3 Z^2 + 1) + \varepsilon^4 X^2 \\
\varepsilon' &= 0
\end{rcases} \label{eq:syscytran12} \noeqref{eq:syscytran12}
\end{align}
where $'$ denotes differentiation with respect to $\bar{\tau}$.
To investigate the degenerate critical manifold branch $\Gamma_1 = \{X=0\}$,  \cref{fig:anasysxz}, we perform a cylindrical blow-up, as shown in \cref{fig:busc1_12}, using the transformation:
\begin{equation}
X = \bar{r} \bar{X}, \quad Z = \bar{Z}, \quad \varepsilon = \bar{r} \bar{\varepsilon}\,,
\end{equation}
under the constraint $\bar{X}^2+\bar{\varepsilon}^2=1$. This cylindrical blow-up unfolds the dynamics near $\Gamma_1$ by mapping the singular origin to the surface of a cylinder. To analyse the vector field on this manifold, we study the system in a series of overlapping \textit{coordinate charts} that cover the different physical regimes of interest.

\subsubsection{Chart $K_1$}

We define the entry chart $K_1$ (see \cref{fig:busc1_12}) by setting $\bar{X} = 1$ and introducing local coordinates $\bar{r} := r_1, \bar{Z} := Z_1, \bar{\varepsilon} := \varepsilon_1$.
This yields the transformation $X = r_1, Z = Z_1, \varepsilon = \varepsilon_1 r_1$ and the desingularised system:
\begin{align}
&\begin{rcases}
r_1' &= - 2  r_1 \varepsilon_1 Z_1\left[ (\nu - \alpha \sigma_1) Z_1^2 - \varepsilon_1^2 \alpha (\sigma_2 r_1^2 + \sigma_3 Z_1^2 + 1) + \varepsilon_1^4 \nu r_1^4 \right]\\
Z_1' &= (1 - \gamma \sigma_1 Z_1) Z_1^2 - \varepsilon_1^2 \gamma Z_1 (\sigma_2 r_1^2 + \sigma_3 Z_1^2 + 1) + \varepsilon_1^4 r_1^4 \\
\varepsilon_1' &= 2 \varepsilon_1^2 Z_1  \left[ (\nu - \alpha \sigma_1) Z_1^2 - \varepsilon_1^2 \alpha (\sigma_2 r_1^2 + \sigma_3 Z_1^2 + 1) + \varepsilon_1^4 \nu r_1^4 \right]
\end{rcases} \label{eq:syscha1} \noeqref{eq:syscha1}
\end{align}
where $'$ denotes differentiation with respect to $\tau_1$ given by $d\bar{\tau} = r_1^{-2} d\tau_1$. Chart $K_1$ covers the transition between Regimes 1 and 2. The transition dynamics are organised by the invariant planes $\varepsilon_1 = 0$ and $r_1 = 0$, which connect at the `equator'.\footnote{Here, the `equator' refers specifically to the two points $\bar{X} = \pm 1$ on the circle $\bar{X}^2 + \bar{\varepsilon}^2 = 1$ (where $\bar{\varepsilon}=0$), extended along the $Z$-axis in the blown-up space.} 

\begin{figure}
\centering
\includegraphics[width=0.85\textwidth]{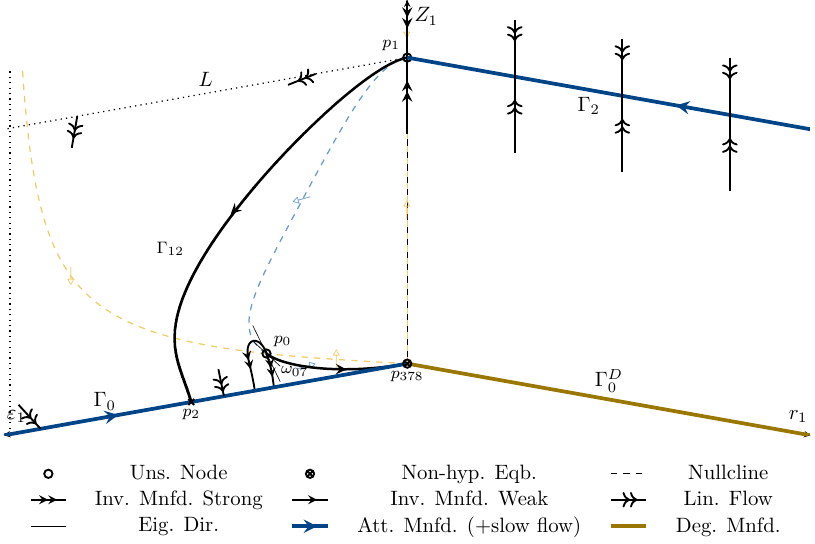}
\caption{Chart $K_1$ with slow-flow overlays.}
\label{fig:3dchartk1}
\end{figure}

\noindent\textbf{Dynamics in invariant plane $\varepsilon_1 = 0$:}
Restricting system \eqref{eq:syscha1} to $\varepsilon_1 = 0$ yields the 2D problem { \eqref{eq:sys:e1cha1}.}
The set of equilibria consists of two lines: $\Gamma_0^D := \{ Z_1 = 0 \}$ and $\Gamma_2 := \{ Z_1 = 1/(\gamma \sigma_1) \}$. 
Importantly, we have gained normal hyperbolicity of $\Gamma_2$ up to $r_1=0$ (point $p_1$ in \cref{fig:3dchartk1}), which lies on the `equator' of the blown-up cylinder; see \cref{app:k1}.
\\

\noindent\textbf{Dynamics in invariant plane $r_1 = 0$:}
Restricting system \eqref{eq:syscha1} to $r_1 = 0$ yields the 2D problem { \eqref{eq:sys:r1cha1} that governs} the transition onto the blown-up cylinder.
The set of equilibria comprises the line $\Gamma_0 := \{Z_1 =0\}$ and two isolated points, $p_0$ and $p_1$, given in $(\varepsilon_1, Z_1)$ coordinates by:
\begin{equation}
p_0 = \left(\sqrt{\frac{\alpha (\nu - \alpha \sigma_1)}{\sigma_3 \alpha^2+ \gamma^2 \nu^2}}, \frac{\alpha}{\gamma \nu} \right)
, \quad
p_1 = \left(0, \frac{1}{\gamma \sigma_1} \right)\,.
\label{eq:ck1isopts}
\end{equation}

The line $\Gamma_0$ (the $\varepsilon_1$-axis) has a non-trivial eigenvalue $\lambda=-\varepsilon_1^2\gamma\le 0$. Hence the line $\Gamma_0$ is normally attracting for $\varepsilon_1>0$ and loses normal hyperbolicity at the origin of the chart, denoted $p_{378} := (0,0)$.

The first equilibrium $p_0$ corresponds to the unstable node identified in Regime 2 (see \cref{sec:regime2} for details). The second equilibrium $p_1$ corresponds to the intersection of the invariant manifold $\Gamma_2$ with the $r_1=0$ plane mentioned above. { Its eigenvalues are given in \cref{app:k1}.}
The local centre manifold $W^c_{loc}(p_1) \subseteq\Gamma_{12}$ is given to leading order by the quadratic approximation { \eqref{app:eq:cmgraph}.}
    The reduced dynamics on this 1D manifold for $\varepsilon_1 \ge 0$ are governed by { \eqref{app:eq:cmflow}.}
    Given $\nu - \alpha\sigma_1 > 0$, the coefficient is positive, indicating that the flow on the manifold is repelling for $\varepsilon_1 > 0$, i.e., moving away from $p_1$, which lies on the `equator', onto the cylinder, see \cref{fig:3dchartk1}.
\begin{proposition}[Existence of the 2D centre manifold]
 The 3D system \eqref{eq:syscha1} possesses a local 2-dimensional centre manifold $W^c_{loc}(p_1)$ at $p_1 = (0, \frac{1}{\gamma \sigma_1}, 0)$, tangent to the $(r_1, \varepsilon_1)$-eigenspace. To leading order, the manifold is represented by the graph:
\begin{equation}
    Z_1 = h(r_1, \varepsilon_1) = \frac{1}{\gamma \sigma_1} - \frac{\sigma_3 + \gamma^2 \sigma_1^2}{\gamma \sigma_1^2} \varepsilon_1^2 + \mathcal{O}(|(r_1, \varepsilon_1)|^3).
\end{equation}

The dynamics on $W^c_{loc}(p_1)$ are organised by two invariant 1D sub-manifolds:
\begin{itemize}
    \item In the invariant plane $\varepsilon_1 = 0$, the normally attracting line of equilibria $\Gamma_2 := \{Z_1 = 1/(\gamma \sigma_1)\}$ guides the Regime 1 slow flow (\cref{sec:regime1reduced}) towards the equator ($r_1 \to 0$).
    \item In the invariant plane $r_1 = 0$, the flow satisfies $\varepsilon_1' > 0$ for $\varepsilon_1 > 0$, ejecting trajectories onto the blown-up cylindrical surface.
\end{itemize}
\end{proposition}
Thus, $p_1$ serves as a transition point: trajectories are attracted towards the neighbourhood of $p_1$ along the $r_1$-direction and are subsequently ejected along the $\varepsilon_1$-direction. {$W^c_{loc}(p_1)$ is pinched within a narrowing wedge between the tangent line $L := \{Z_1 = 1/(\gamma\sigma_1)\}$ and the $Z_1$-nullcline, along which $Z_1$ decreases; see \cref{app:k1}.} The fate of this centre manifold is established in chart $K_2$.

\subsubsection{Chart $K_2$}

We next analyse the dynamics in the scaling chart $K_2$ (\cref{fig:busc1_12}), defined by setting $\bar{\varepsilon} = 1$ with local coordinates $(X_2, Z_2, r_2)$. This yields the transformation $X = r_2 X_2, Z = Z_2, \varepsilon = r_2$ and, after desingularisation, the system:
\begin{align}
&\begin{rcases}
X_2' &= -2 Z_2 [(\nu - \alpha \sigma_1) X_2^2 Z_2^2 - \alpha (\sigma_3 Z_2^2 + 1) - r_2^2 \alpha \sigma_2 X_2^2 + r_2^4 \nu X_2^2] \\
Z_2' &= (1 - \gamma \sigma_1 Z_2) X_2^2 Z_2^2 - \gamma Z_2 (\sigma_3 Z_2^2 + 1) - r_2^2 \gamma \sigma_2 X_2^2 Z_2 + r_2^4 X_2^2 \\
r_2' &= 0
\end{rcases} \label{eq:syscha2}
\end{align}
where $'$ denotes differentiation with respect to $\tau_2$ ($d\bar{\tau} = r_2^{-2} d\tau_2$). Chart $K_2$ fully resolves Regime 2, including the $Z$-axis ($X_2=0$) which is inaccessible in Chart $K_1$. Notably, in the singular limit $r_2=0$, system \eqref{eq:syscha2} is identical to the System $U$-$Z$ \eqref{eq:spuz}, which has already been analysed in \cref{sec:regime2}. 

The results from the $U$-$Z$ analysis -- specifically the existence of the unstable node $p_0$ and the attracting critical manifold $\Gamma_5 := \{Z_2=0\}$ -- directly determine the global fate of the centre manifold $W^c(p_1)$ arriving from the entry chart $K_1$. While the local analysis in $K_1$ shows the trajectory is trapped near the equator, the global topology of the layer problem in $K_2$ confirms that the flow is repelled by $p_0$ and subsequently turned by the $X_2$-nullcline. This forces a fast heteroclinic connection on the blown-up cylinder from $p_1$ towards a point $p_2$ on the critical manifold $\Gamma_5$; see \cref{fig:anasysuz} and \cref{fig:3dchartk1}.

\subsubsection{Summary of the First Blow-up}

The dynamics following the cylindrical blow-up of $\Gamma_1 = \{X=0\}$ are obtained by superimposing charts $K_1$ and $K_2$, effectively desingularising the transition from Regime 1 to Regime 2. Chart $K_1$ establishes the entry of the flow from the normally hyperbolic attracting manifold $\Gamma_2$. The resulting centre manifold $W^c_{loc}(p_1)$ acts as the primary gateway from the macroscopic state space into the singular interior.

As revealed by the global topology in chart $K_2$, the connection from $p_1$ to $p_2$ is enforced by a geometric trapping region: the flow is pinched between the tangent line $L$ and the $Z_1$-nullcline (blue) locally, while being repelled by the unstable node $p_0$ globally. This sequence establishes a heteroclinic connection on the blown-up cylinder, mapping the trajectory from the equator onto the critical manifold $\Gamma_5$. Subsequently, the slow flow $r_2$ on $\Gamma_5$ (corresponding to the Regime 2 degenerate reduced problem) carries the trajectory towards the equilibrium $p_{378}$ (see \cref{fig:3dchartk1}), completing the transition across this first degenerate region and setting the stage for the second blow-up.

\subsection{Second cylindrical blow-up}

\begin{figure}
\centering
\includegraphics[width=0.9\textwidth]{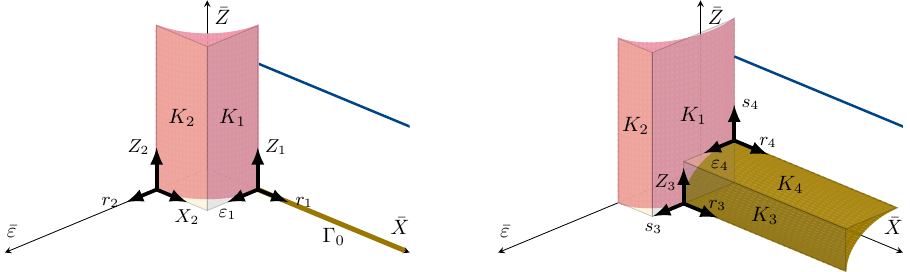}
\caption{Second cylindrical blow-up geometry.}
\label{fig:busc1d_1234}
\end{figure}

To resolve the degenerate dynamics near the non-hyperbolic point $p_{378} = (0,0)$ in chart $K_1$ and describe the transition into the $Z \sim \varepsilon^2$ scaling of Regime 3, we perform a second, weighted cylindrical blow-up. We introduce the transformation:
\begin{equation}
r_1 = \hat{r}_1, \quad Z_1 = \hat{s}^2 \hat{Z}_1, \quad \varepsilon_1 = \hat{s} \hat{\varepsilon}_1\,,
\end{equation}
subject to the constraint $\hat{Z}_1^2+\hat{\varepsilon}_1^2=1$. This quasi-homogeneous scaling desingularises the vector field at $p_{378}$, unfolding the singular origin into a cylindrical surface that connects the Regime 2 transition to the Regime 3 exit. The geometry of this second blow-up is presented in \cref{fig:busc1d_1234}, and its dynamics are analysed via the transition chart $K_3$ and the exit chart $K_4$.

\subsubsection{Chart $K_3$}

The transition chart $K_3$ (see \cref{fig:busc1d_1234}) is defined by setting $\hat{\varepsilon}_1 = 1$, with local coordinates $(r_3, Z_3, s_3)$. This yields the quasi-homogeneous transformation $r_1 = r_3, Z_1 = s_3^2 Z_3, \varepsilon_1 = s_3$, and the desingularised system:
\begin{align}
&\begin{rcases}
r_3' &= -2 Z_3 r_3 s_3^3 \left[s_3^2 (\nu -\alpha (\sigma_1 + s_3^2 \sigma_3) Z_3^2 + \nu r_3^4) -\alpha (\sigma_2 r_3^2 + 1) \right] \\
Z_3' &= Z_3^2 + r_3^4 - \gamma (\sigma_2 r_3^2 + 1) Z_3 - s_3^2 \gamma Z_3^3 (\sigma_1 + s_3^2 \sigma_3) \\
&- 4 s_3^3 Z_3^2 \left[s_3^2 (\nu - \alpha (\sigma_1 + s_3^2 \sigma_3) Z_3^2 +\nu r_3^4) -\alpha (\sigma_2 r_3^2 + 1) \right] \\
s_3' &= 2 Z_3 s_3^4 \left[s_3^2 (\nu - \alpha (\sigma_1 + s_3^2 \sigma_3) Z_3^2 +\nu r_3^4) -\alpha (\sigma_2 r_3^2 + 1) \right]
\end{rcases} \label{eq:syscha3} \noeqref{eq:syscha3}
\end{align}
where $'$ denotes differentiation with respect to $\tau_3$ ($d\bar{\tau} = s_3^{-2} d\tau_3$). Chart $K_3$ resolves the transition from Regime 2 to Regime 3. The dynamics on the blown-up cylinder are organised by the invariant planes $s_3 = 0$ (the surface of the cylinder) and $r_3 = 0$ (the blown-up $Z$-axis).

\noindent\textbf{Dynamics in invariant plane $s_3 = 0$:}
Restricting the system to the invariant plane $s_3=0$ yields a 2D flow where $r_3$ acts as a constant of motion ($r_3' = 0$). The vertical dynamics on the blown-up cylinder are governed by {\eqref{app:eq:z3quad}.}
The set of equilibria forms the critical manifold $\Gamma_4$, a parabolic curve in the $(r_3, Z_3)$ plane defined by the roots of the quadratic. This manifold is composed of an attracting branch $\Gamma_{4a}$ and a repelling branch $\Gamma_{4r}$, which meet at the fold point $p_4$:
\begin{equation}
p_4 = \left(\sqrt{\frac{\gamma}{2 - \gamma \sigma_2}}, \frac{\gamma}{2 - \gamma \sigma_2}\right),
\end{equation}
{where the eigenvalues vanish; see \cref{app:k3}.}
Importantly, the attracting branch $\Gamma_{4a}$ extends to the point $p_3$ at the origin.\\

\noindent\textbf{Dynamics in invariant plane $r_3 = 0$:}
Restricting the system to the invariant plane $r_3=0$ reveals the transition from the internal dynamics of the first blow-up to the fold geometry. The 2D flow is given by {\eqref{app:eq:r3cha3}.}
Factoring out $Z_3$ shows that the set of equilibria comprises the line $\Gamma_5 := \{Z_3=0\}$ (the continuation of the attracting manifold from Regime 2) and the isolated point $p_7 = (0, \gamma)$ in the $(s_3, Z_3)$ plane.

The point $p_7 = (0, \gamma)$ is a partially hyperbolic equilibrium. {Its linearisation and centre flow are given in \cref{app:k3}.} For $s_3 > 0$, the centre flow is attracted towards $p_7$ ($s_3' < 0$). This configuration establishes a heteroclinic connection from the unstable node $p_0$ (inherited from chart $K_1$, whose weak unstable manifold on the right connects to the non-hyperbolic point $p_{378}$) to $p_7$. This connection, $\omega_{07}$, serves as a link to the repelling manifold $\Gamma_{4r}$.

The manifold $\Gamma_5$ is normally hyperbolic and attracting with eigenvalue $\lambda = -\gamma < 0$ up to the point $p_3$ at the origin, successfully resolving the loss of hyperbolicity at $p_{378}$ in the original coordinates. The point $p_3 = (0, 0, 0)$ in the full 3D system is a non-isolated equilibrium point representing the corner where the first blow-up scaling ends and the second begins.

\begin{proposition}[Existence of the 2D centre manifold at $p_3$]
The 3D system \eqref{eq:syscha3} possesses a local 2-dimensional centre manifold $W^c_{loc}(p_3)$ at $p_3 = (0, 0, 0)$, tangent to the $(r_3, s_3)$-eigenspace. To leading order, the manifold is represented by the graph $Z_3 = \Phi(r_3, s_3) = \mathcal{O}(r_3^4)$.
The dynamics on $W^c_{loc}(p_3)$ are organised by the two invariant planes that intersect at $p_3$:
\begin{itemize}
    \item In the invariant plane $r_3 = 0$, the manifold contains the equilibrium branch $\Gamma_5 := \{Z_3 = 0\}$. The flow along $\Gamma_5$ is directed \textit{towards} $p_3$.
    \item In the invariant plane $s_3 = 0$, the manifold contains the attracting branch $\Gamma_{4a}$. The flow on this branch is directed \textit{away} from $p_3$.
\end{itemize}
\end{proposition}

\begin{figure}
\centering
\includegraphics[width=0.85\textwidth]{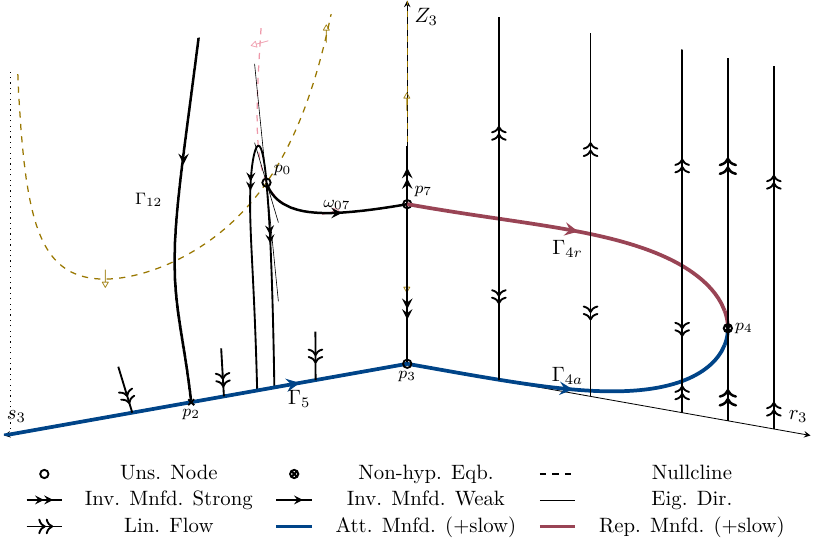}
\caption{Chart $K_3$ with slow-flow overlays.}
\label{fig:3dchartk3}
\end{figure}

\subsubsection{Chart $K_4$}

The exit chart $K_4$ (see \cref{fig:busc1d_1234}) is defined by setting $\hat{Z}_1 = 1$ with local coordinates $(r_4, s_4, \varepsilon_4)$, i.e., the transformation $r_1 = r_4, Z_1 = s_4^2, \varepsilon_1 = s_4 \varepsilon_4$ and, after desingularisation, the system:
\begin{align}
&\begin{rcases}
r_4' &= - 4 r_4 s_4^3 \varepsilon_4 \left[ (\nu - \alpha \sigma_1) s_4^2 - \varepsilon_4^2 \alpha (\sigma_2 r_4^2 + \sigma_3 s_4^4 + 1 ) + \varepsilon_4^4 \nu r_4^4 s_4^2 \right]\\
s_4' &= s_4 \left[(1 - \gamma \sigma_1 s_4^2) - \varepsilon_4^2 \gamma (\sigma_2 r_4^2 + \sigma_3 s_4^4 + 1) + \varepsilon_4^4 r_4^4 \right] \\
\varepsilon_4' &= - \varepsilon_4 [(1 - \gamma \sigma_1 s_4^2) - \varepsilon_4^2 \gamma (\sigma_2 r_4^2 + \sigma_3 s_4^4 + 1) + \varepsilon_4^4 r_4^4 \\
&- 4 s_4^3 \varepsilon_4 ( (\nu - \alpha \sigma_1) s_4^2 - \varepsilon_4^2 \alpha (\sigma_2 r_4^2 + \sigma_3 s_4^4 + 1 ) + \varepsilon_4^4 \nu r_4^4 s_4^2 ) ]
\end{rcases} \label{eq:syscha4} \noeqref{eq:syscha4}
\end{align}
where $'$ denotes differentiation with respect to $\tau_4$ ($d\bar{\tau} = 2 s_4^{-2} d\tau_4$). Chart $K_4$ resolves the transition from the singular scaling of Regime 3 back to the macroscopic dynamics of Regime 1. The dynamics on the blown-up cylinder are organised by the invariant planes $r_4=0$ (blown-up $Z$-axis), $s_4=0$ (surface of the cylinder), and $\varepsilon_4=0$ (Regime 1).

\noindent\textbf{Dynamics in invariant plane $r_4 = 0$:}
Restricting the system to the invariant plane $r_4 = 0$ yields the 2D flow {\eqref{app:eq:r4cha4}.}
The equilibria in this $(\varepsilon_4, s_4)$-plane include the unstable node $p_0$, the saddle-node $p_1$, and the gateway point $p_7 = (1/\sqrt{\gamma}, 0)$, analysed in the previous chart. The only new equilibrium point introduced in this chart is the origin $p_8 = (0,0)$, which is a hyperbolic saddle{; see \cref{app:k4} for details.} This saddle structure acts as a partition, directing trajectories away from the singular scaling and back towards the macroscopic Regime 1 dynamics as $s_4$ increases.

The combination of charts $K_3$ and $K_4$ reveals that the loss of normal hyperbolicity associated with the singularity $p_{378}$ is fully resolved. This degenerate point unfolds into a sequence of three desingularised structures: the normally hyperbolic manifold $\Gamma_5$ (in $K_3$), the partially hyperbolic gateway $p_7$, and the hyperbolic saddle $p_8$ (in $K_4$). 

\noindent\textbf{Dynamics in invariant plane $s_4 = 0$:}
Restricting the system to the invariant plane $s_4 = 0$ (the surface of the cylinder in the exit chart) yields a 2D flow where $r_4$ is a constant of motion. The vertical dynamics are governed by { \eqref{app:eq:s4cha4}.}
The set of equilibria includes the normally attracting line $\Gamma_0^S := \{ \varepsilon_4 = 0 \}$ and the folded critical manifold $\Gamma_4$, defined by the roots of the quadratic {\eqref{app:eq:g4quart}} in $\varepsilon_4^2$.
This confirms that the manifold $\Gamma_4$ identified in chart $K_3$ persists in the exit chart. In this representation, $p_7$ is recovered at $(r_4, \varepsilon_4) = (0, 1/\sqrt{\gamma})$. The fold point $p_4$, where the manifold branches meet and normal hyperbolicity is lost, is located at:
\begin{equation}
p_4 = \left(\sqrt{\frac{\gamma}{2-\gamma \sigma_2}}, \sqrt{\frac{2- \gamma \sigma_2}{\gamma}}\right).
\end{equation}
The persistence of this parabolic structure across both $K_3$ and $K_4$ ensures a continuous transition of the manifold geometry, rigorously establishing the infra-slow crawl and subsequent jump that characterise Regime 3.

\begin{figure}
\centering
\includegraphics[width=0.75\textwidth]{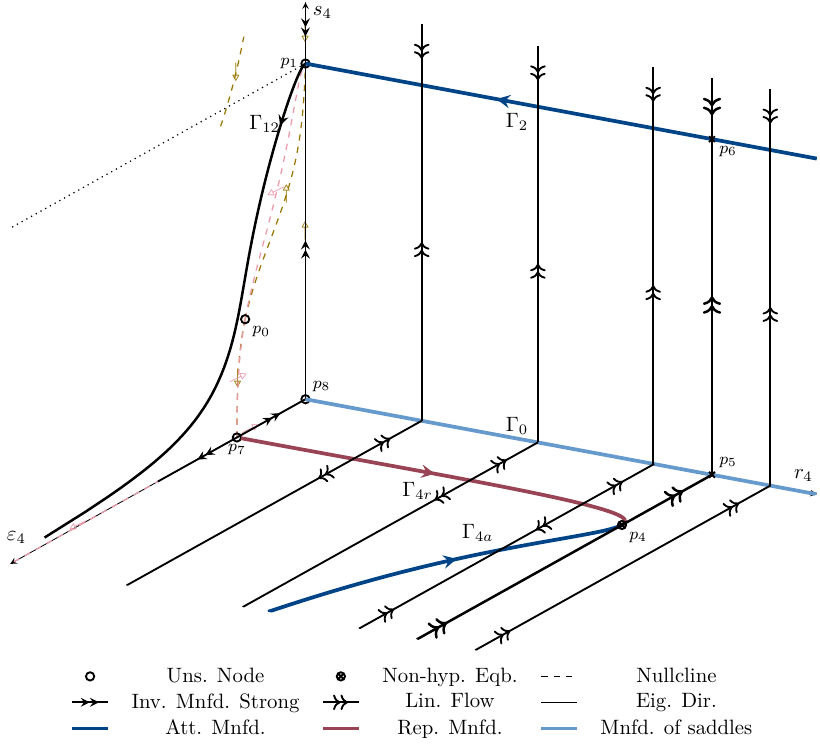}
\caption{Chart $K_4$ with slow-flow overlays.}
\label{fig:3dchartk4}
\end{figure}

\noindent\textbf{Dynamics in invariant plane $\varepsilon_4 = 0$:}
Restricting the system to $\varepsilon_4 = 0$ yields $r_4'=0$ and the 1D flow:
\begin{equation}
s_4' = s_4 (1 - \gamma \sigma_1 s_4^2).
\end{equation}
This plane recovers the primary Regime 1 skeleton, serving as the `landing zone' for trajectories exiting the desingularised transition region. The invariant line $s_4 = 0$ captures the local trace of the blown-up repelling critical manifold $\Gamma_{4r}$ (the $X$-axis in original coordinates, see \cref{fig:anasysuv}). More significantly, we recover the stable manifold $\Gamma_2$ at $s_4 = 1/\sqrt{\gamma \sigma_1}$.\footnote{The square-root relationship between the coordinate $s_4 = 1/\sqrt{\gamma \sigma_1}$ and the $K_1$ equilibrium $Z_1 = 1/(\gamma \sigma_1)$ is a consequence of the weighted transformation $Z_1 = s_4^2$ used in the second blow-up.} 
As illustrated in \cref{fig:3dchartk4}, this recovery establishes the global consistency of the analysis: trajectory segments jumping from the fold point are successfully reinjected into the neighbourhood of the attracting manifold $\Gamma_2$, completing the heteroclinic loop and returning the flow to the macroscopic dynamics of Regime 1.

\noindent\textbf{The 3D fast jump and Dulac transition:}
While the invariant planes form the skeleton of the exit chart, completing the heteroclinic loop requires the full 3D flow. Upon reaching the fold $p_4$ in the $s_4=0$ plane, normal hyperbolicity is lost and the singular trajectory is discharged along the fast fibre (decreasing $\varepsilon_4$) until it lands on the equatorial point:
\begin{equation}
p_5 = \left(\sqrt{\frac{\gamma}{2-\gamma \sigma_2}}, 0, 0\right).
\end{equation}
This point lies at the intersection of the invariant planes on the equator, which forms a normally hyperbolic line of equilibria. Linearising the full 3D system \eqref{eq:syscha4} at $p_5$ gives the eigenvalues $(0, 1, -1)$. The stable fibre associated with the eigenvalue $-1$ (lying in the $s_4=0$ plane) attracts the singular flow from $p_4$, while the unstable fibre associated with the eigenvalue $+1$ (lying in the $\varepsilon_4=0$ plane) is responsible for ejecting the flow towards $p_6$ on $\Gamma_2$. Since the (absolute) ratio of these non-trivial eigenvalues is $1$ (a $1:1$ resonance), the passage near this equatorial point $p_5$ generates a Dulac transition map with logarithmic corrections of order $\mathcal{O}(\varepsilon \ln(1/\varepsilon))$. As illustrated in \cref{fig:3dchartk4}, this transition successfully reinjects the trajectory into the macroscopic dynamics of Regime 1, and the resulting Dulac map drives the persistence analysis in \cref{prop:persistence}.

\subsubsection{Summary of the Second Blow-up}

The analysis of charts $K_3$ and $K_4$ successfully desingularises the flow through the degenerate point $p_{378}$, unfolding this singularity into a sequence of hyperbolic and partially hyperbolic structures. The local flow tracking the singular limit cycle within this blow-up is organised as follows:
\begin{itemize}
    \item \textit{Entry ($p_3$):} The manifold $\Gamma_5$ (accessible in chart $K_3$) provides the entry path from Regime 2. The flow along $\Gamma_5$ carries trajectories directly into the non-isolated equilibrium $p_3$ at the origin of the chart.
    \item \textit{Manifold connection ($p_3 \to p_4$):} The local centre manifold at $p_3$ geometrically connects the incoming flow from $\Gamma_5$ to the attracting branch of the folded manifold $\Gamma_{4a}$. The trajectory moves away from $p_3$ along $\Gamma_{4a}$ and accelerates towards the fold point $p_4$.
    \item \textit{Exit mechanism (fast jump):} Upon reaching the fold point $p_4$, normal hyperbolicity is lost and the trajectory is discharged into the fast layer flow. Chart $K_4$ resolves the full 3D topology of this exit: the trajectory drops to the normally hyperbolic equatorial point $p_5$, where the unstable fibre ejects it out of the singular scaling and back towards the macroscopic attracting manifold $\Gamma_2$ (at $p_6$). 
\end{itemize}
This local desingularisation provides the mathematical foundation for the global singular limit cycle. Furthermore, the $1:1$ eigenvalue resonance at $p_5$ generates the Dulac transition map that, combined with the fold asymptotics, dictates the bounds in the final proof of existence for the periodic orbit.
\subsection{The singular limit cycle}

\begin{figure}
\centering
\includegraphics[width=0.9\textwidth]{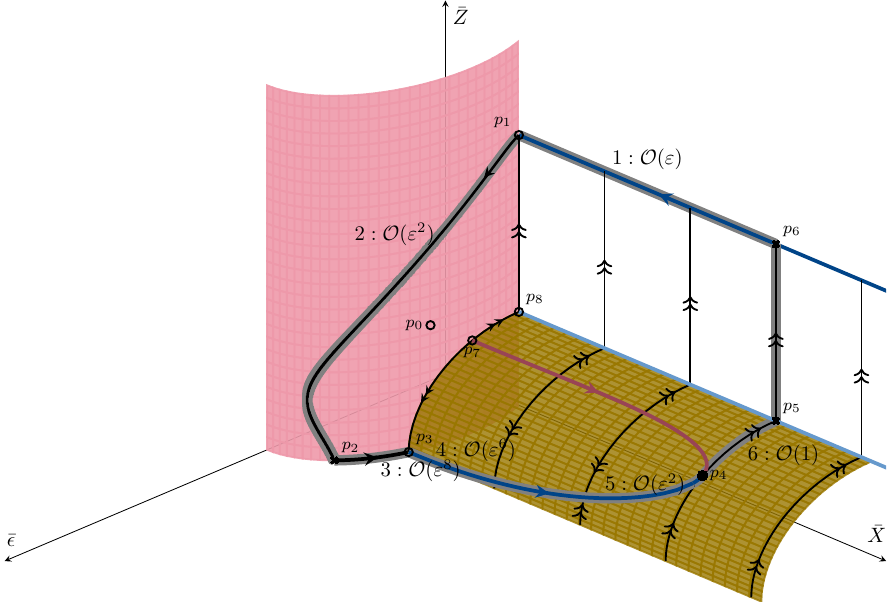}
\caption{Singular limit cycle. A 3D schematic of the full singular limit cycle in the blown-up space, showing the key points ($p_1, p_4, p_6$, etc.) and the path of the trajectory as it moves between the slow manifolds and fast fibers.}
\label{fig:3dbusl}
\end{figure}

The full singular limit cycle $\Gamma^0$ is a composite trajectory constructed from the segments analysed in the GSPT `patchwork' and geometrically `glued' together by the blow-up analysis. As illustrated in \cref{fig:3dbusl}, the cycle is rigorously traced through the sequence of connections $p_1 \to p_2 \to p_3 \to p_4 \to p_5 \to p_6 \to p_1$.
This path connects the various critical manifolds and is defined by a distinct hierarchy of timescales that explains the oscillator's physical phases:

\begin{enumerate}
    \item \textit{Slow drift ($\mathcal{O}(\varepsilon)$):} Starting at $p_6$ on $\Gamma_2$, the trajectory moves towards the $Z$-axis.
    \item \textit{Slower deceleration ($\mathcal{O}(\varepsilon^2)$):} At $p_1$, the flow is captured by the first blow-up, guiding the trajectory down the $Z$-axis.
    \item \textit{Slowest crawl ($\mathcal{O}(\varepsilon^8)$):} Near the origin, the flow enters the stagnation phase. This represents the long `refill' phase of the two-stroke cycle.
    \item \textit{Infra-slow acceleration ($\mathcal{O}(\varepsilon^6)$):} Transitioning onto $\Gamma_{4a}$, the trajectory begins to accelerate as it moves along it.
    \item \textit{Slower acceleration ($\mathcal{O}(\varepsilon^2)$):} Approaching the fold $p_4$, the Regime 3 layer problem dominates, and the trajectory gains significant speed.
    \item \textit{Fast jump ($\mathcal{O}(1)$):} At $p_4$, the Regime 1 layer problem takes over, causing a near-instantaneous jump past the equatorial saddle $p_5$ and resetting the trajectory back to $p_6$ on $\Gamma_2$.
\end{enumerate}

\subsection{Persistence of the periodic orbit}

Having constructed the singular cycle $\Gamma^0$, we now establish its persistence for sufficiently small $\varepsilon > 0$. Specifically, we will show that $\Gamma^0$ perturbs to a nearby, genuine periodic orbit $\Gamma_\varepsilon$. 

The classical route to establishing persistence relies on continuation from $\varepsilon=0$. The implicit function theorem, applied to the fixed-point equation $P_\varepsilon(h)-h=0$, requires an unperturbed map $P_0$ with an invertible linearisation. This approach fails here. First, no unperturbed map exists: transit times along the slow manifolds diverge as $\varepsilon \to 0^+$. Second, while the return map $P_\varepsilon$ and its passage estimates are uniform for small $\varepsilon$, its dependence on $\varepsilon$ is not smooth at $\varepsilon=0$. The passage past the fold introduces a fractional term of $\mathcal{O}(\varepsilon^{2/3})$, and the passage past the equatorial saddle introduces logarithmic terms $\mathcal{O}(\ln(1/\varepsilon))$, whose $\varepsilon$-derivatives diverge as $\varepsilon\to 0^+$.

Although working in the $\varepsilon$-extended phase space and employing sequential blow-up transformations (\cref{sec:blowup}) desingularises the flow geometry, the loss of parameter differentiability inherent to these transitions prevents the application of the implicit function theorem with respect to the original $\varepsilon$ coordinate.

Since smoothness in $\varepsilon$ is violated, the family must be constructed pointwise. For each fixed $\varepsilon \in (0, \varepsilon_0]$, we define a smooth cross-section $\Sigma$ and a return map $P_\varepsilon(h)$, where $h$ measures the transverse distance to $\Gamma^0$. The blow-up analysis provides contraction estimates that are uniform in $\varepsilon$. This allows us to prove that $P_\varepsilon$ is a strict contraction on a parameter-independent neighbourhood $I$, supplying the fixed point via the contraction mapping theorem.

\begin{proposition}[Persistence]\label{prop:persistence}
There exist $\varepsilon_0>0$ and a family $\{\Gamma_\varepsilon\}_{\varepsilon\in(0,\varepsilon_0]}$ of periodic orbits, one for each $\varepsilon$, such that the Hausdorff distance 
\[ 
\dist_{\mathrm H}\big(\Gamma_\varepsilon,\Gamma^0\big)=\mathcal{O}(\varepsilon^{2/3})\longrightarrow 0 \qquad\text{as }\varepsilon\to0^+ . 
\] 
Each $\Gamma_\varepsilon$ is the unique periodic orbit in a fixed neighbourhood of $\Gamma^0$, and is hyperbolic and exponentially stable, with non-trivial Floquet multiplier $0<\mu_\varepsilon\le e^{-C/\varepsilon}$ for some $C>0$.
\end{proposition}

\begin{proof}
Fix $\varepsilon\in(0,\varepsilon_0]$. As established in \cref{sec:returnmap}, the return map $P_\varepsilon$ is a self-map on a compact, parameter-independent neighbourhood $I$ of $\Gamma^0\cap\Sigma$. Because the local passage bounds are uniform in $\varepsilon$, the overall linearisation satisfies $|\partial_h P_\varepsilon|\le e^{-C/\varepsilon}<1$ for all $h \in I$. 

$P_\varepsilon$ is therefore a strict contraction on the complete metric space $I$. By the contraction mapping theorem, it has a unique fixed point $h^*_\varepsilon\in I$; the orbit through $h^*_\varepsilon$ returns to itself after one transit and is therefore closed. This yields the periodic orbit $\Gamma_\varepsilon$, and its uniqueness in $I$ establishes uniqueness in a fixed neighbourhood of $\Gamma^0$. 

Repeating this for each $\varepsilon \in (0, \varepsilon_0]$ produces the family. The exponential stability and the Floquet multiplier bound follow directly from the strong contraction property $|\partial_h P_\varepsilon| \le e^{-C/\varepsilon} < 1$. Finally, the convergence rate $\dist_{\mathrm H}(\Gamma_\varepsilon,\Gamma^0) = \mathcal{O}(\varepsilon^{2/3})$ follows from the maximal transverse excursion occurring at the fold passage $p_4$, as derived in \cref{sec:returnmap}.
\end{proof}

\section{Conclusion}
\label{sec:discussion}

The primary objective of this work was to rigorously analyse the dynamics of the Smolen model for glycolytic oscillations, a key component of the Integrated Oscillator Model for the activity of pancreatic \bcells.  The analysis presented here overcomes three fundamental mathematical hurdles that previously obscured the multi-scale nature of this system.

The first challenge was the \textit{non-dimensionalisation} of the original biochemical rate equations. By identifying the relative ratios of maximal enzymatic velocities and substrate affinities, we were able to define a key small parameter, $\varepsilon$, which serves as the ratio of disparate metabolic timescales. This step was essential to recast the glycolytic oscillation as a singular perturbation problem.

Building on this, we addressed the fact that the biophysical model \eqref{eq:rndxy} is not directly amenable to standard GSPT analysis due to its non-polynomial form and the singularity at $Y=0$. We applied a state-dependent time rescaling to transform the system into a polynomial `surrogate' relaxation oscillator \eqref{eq:NdSysXZ}. The GSPT and blow-up analyses were thus applied to this surrogate system as a necessary method to resolve the singular geometry that is hidden in the original model.

The analysis of the surrogate system revealed a vast hierarchy of timescales. Because the time rescaling is state-dependent, these surrogate timescales -- stretching from $\mathcal{O}(1)$ to $\mathcal{O}(\varepsilon^8)$ -- do not correspond directly to speeds in the biophysical system. However, they provide the essential topological map required to construct a singular limit cycle $\Gamma^0$. Our results show that the glycolytic oscillation is not a simple binary switch, but a multi-layered hierarchy of regimes.

This work significantly extends the analysis of glycolytic oscillators.
In particular, the dynamics of the Smolen model were analysed within the context of the IOM in \cite{McKenna2018}, but the wide range of timescales was not identified. While our methods were inspired by foundational work on the Goldbeter-Lefever model for glycolytic oscillations in yeast \cite{Kosiuk2011}, the glycolytic oscillations in \bcells\ are due to a different biochemical mechanism: positive feedback and substrate depletion. We show that this model generates a far more intricate singular geometry and a richer timescale hierarchy than the models previously studied. The `pulse' segment of the glycolytic oscillation consists of flow through three regimes in the analysis, while the `refill' segment is driven by ultra-slow/slow flow through Regime 3. This provides a biophysical example of a two-stroke relaxation oscillation \cite{LeCorbeiller1960}. 

\backmatter

\bmhead{Acknowledgements}
The authors would like to thank Patrick Fletcher for many insightful
discussions.

\section*{Statements and Declarations}

\bmhead{Competing interests}
The authors declare no competing interests.

\bmhead{Funding}
This work was supported by the Australian Research Council Discovery Project DP220101817.

\bmhead{Data availability}
No datasets were generated or analysed during the current study. All results are derived analytically or obtained by numerical solution of the model equations; the equations and parameter values needed to reproduce every figure are given in the article (Appendix~\ref{app:sec:bio}).

\bmhead{Author contributions}
\textbf{Prannath Moolchand:} Conceptualization (equal); Methodology (lead); Formal analysis (lead); Software (lead); Investigation (lead); Visualization (lead); Writing -- original draft (lead); Writing -- review and editing (equal).
\textbf{Richard Bertram:} Conceptualization (supporting); Writing -- review and editing (supporting).
\textbf{Martin Wechselberger:} Conceptualization (equal); Methodology (supporting); Formal analysis (supporting); Supervision (lead); Funding acquisition (lead); Writing -- review and editing (equal).

\newpage
\begin{appendices}

\crefalias{section}{appendix} 
\crefalias{subsection}{appendix}
\crefalias{subsubsection}{appendix}

\section{Glycolytic flux functions}
\label{app:sec:bio}

The metabolic flux terms in the glycolytic model (in $\mu$M/ms) are:
\begin{align}
J_{GK} &= \dfrac{v_{GK}}{1 + (k_{gk}/g_{lce})^{h_{gkglc}}}, \qquad
J_{PDH} = v_{PDH}\sqrt{\dfrac{\text{FBP}}{1 \mu M}} \\
J_{PFK} &= v_{PFK}\dfrac{(1-k_{PFK})w_{1110} + k_{PFK}\sum_{i,j,l \in \{0,1\}}w_{ij1l}}{\sum_{i,j,k,l \in \{0,1\}}w_{ijkl}} \label{app:eq:jpfk}
\end{align}

The weights $w_{ijkl}$, given by
\begin{equation}
w_{ijkl} = \dfrac{(\text{AMP}/K_1)^i(\text{FBP}/K_2)^j(\text{F6P}^2/K_3^2)^k(\text{ATP}^2/K_4^2)^l}{f_{13}^{ik}f_{23}^{jk}f_{41}^{il}f_{42}^{jl}f_{43}^{kl}}, \label{app:eq:weights} \noeqref{app:eq:weights}
\end{equation}
depend on the metabolites $\text{AMP}$, $\text{ADP}$, and $\text{ATP}$, which are related by:
\begin{equation}
\text{ADP} = \dfrac{\sqrt{\text{ATP}(4 A_{tot}-3 \text{ATP})}-\text{ATP}}{2}, \qquad
\text{AMP} = \dfrac{\text{ADP}^2}{\text{ATP}}. \label{app:eq:amp}
\end{equation}

The biophysical parameters and their units are summarised in \cref{app:tab:physparams,app:tab:sysxy}.

\begin{talltblr}[
  caption = {Parameters for System \eqref{app:eq:f6pfbp}.},
  label   = {app:tab:physparams},
]{
  colspec = {c}
}
\resizebox{\textwidth}{!}{
$\begin{tblr}{
  colspec = {c|l|l||c|l|l||c|l|l},
  columns = {mode=math}
}
\text{Par} & \text{Val} & \text{Unit} & \text{Par} & \text{Val} & \text{Unit} & \text{Par} & \text{Val} & \text{Unit} \\
\hline
K_{GPI}   & \num{3.33}  & 1         & v_{GK}    & \num{1.11e-2}  & \mu M/ms  & K_3       & \num{2.24e2}  & \mu M   \\
K_{LG}    & \num{3e-1}  & 1         & v_{PDH}   & \num{1.04e-3}  & \mu M/ms  & K_4       & \num{3.16e1}  & \mu M   \\
k_{gk}    & \num{1.3e1} & 1         & v_{PFK}   & \num{1.04e-2}  & \mu M/ms  & f_{13}    & \num{2e-2}    & 1       \\
g_{lce}   & \num{7}     & 1         & A_{tot}   & \num{3e3}      & \mu M     & f_{23}    & \num{2e-1}    & 1       \\
h_{gkglc} & \num{4}     & 1         & K_1       & \num{3e1}      & \mu M     & f_{41}    & \num{2e1}     & 1       \\
k_{PFK}   & \num{6e-2}  & 1         & K_2       & \num{1}        & \mu M     & f_{42}    & \num{2e1}     & 1       \\
          &             &           &           &                &           & f_{43}    & \num{2e1}     & 1
\end{tblr}$
}
\end{talltblr}

\begin{talltblr}[
  caption = {Parameters for biophysical System \eqref{eq:sysxy}.},
  label   = {app:tab:sysxy},
]{
  colspec = {c}
}
\resizebox{\textwidth}{!}{
$\begin{tblr}{c|l|l||c|l|l||c|l|l}
\text{Par} & \text{Val} & \text{Unit} & \text{Par} & \text{Val} & \text{Unit} & \text{Par} & \text{Val} & \text{Unit} \\
\hline
\omega    & \num{1.00e+00} & \mu M     & \nu       & \num{1.04e-02} & \mu M/ms  & \kappa_3  & \num{3.61e-03} & \mu M     \\
\beta     & \num{2.31e-01} & 1         & \gamma    & \num{5.20e-04} & \mu M/ms  & \kappa_4  & \num{1.89e-04} & 1         \\
\eta      & \num{7.69e-01} & 1         & \kappa_1  & \num{1.12e+01} & (\mu M)^3 & \kappa_5  & \num{1.55e+01} & (\mu M)^3 \\
\alpha    & \num{8.61e-04} & \mu M/ms  & \kappa_2  & \num{8.51e+00} & (\mu M)^2 & \kappa_6  & \num{1.42e+02} & (\mu M)^2
\end{tblr}$
}
\end{talltblr}

\section{Analysis of the dimensionless model}
\label{app:sec:linana}

The Jacobian of system \eqref{eq:rndxy} is given by 
$J = \begin{pmatrix} - \hat{\nu}_1 \hat{r}_X & - \hat{\nu}_1 \hat{r}_Y \\ \hat{\nu}_2 \hat{r}_X & \hat{\nu}_2 \hat{r}_Y - \hat{\gamma}/(2 \sqrt{Y}) \end{pmatrix}$. 
For $X > 0$, the partial derivatives 
\[
\hat{r}_X = \frac{2 X (Y + \hat{\sigma}_6 )(\hat{\sigma}_3 Y + \hat{\sigma}_4)}{D^2}, \qquad \hat{r}_Y = \frac{X^2 \left[ (\hat{\sigma}_2 - \hat{\sigma}_1 \hat{\sigma}_6 ) X^2 + \hat{\sigma}_4 - \hat{\sigma}_3 \hat{\sigma}_6 \right]}{D^2}
\]
are strictly positive, where $D = \hat{\sigma}_1 X^2 Y + \hat{\sigma}_2 X^2 + \hat{\sigma}_3 Y + \hat{\sigma}_4 > 0$. The positivity of $\hat{r}_Y$ holds due to the parameter constraints $\hat{\sigma}_2 > \hat{\sigma}_1 \hat{\sigma}_6$ and $\hat{\sigma}_4 > \hat{\sigma}_3 \hat{\sigma}_6$. 
Evaluated at the equilibrium point $(X_{eq}, Y_{eq})$ \eqref{eq:eqbmxy}, where $\hat{r}^*_X, \hat{r}^*_Y > 0$, the determinant, trace, and discriminant yield:
\begin{align}
\det(J) &= \frac{\hat{\nu}_1 \hat{\gamma} \hat{r}^*_X}{2 \sqrt{Y_{eq}}} > 0, \qquad 
\tr(J) = \hat{r}^*_Y - \hat{\nu}_1 \left( \hat{r}^*_X + \frac{\hat{\gamma}^2}{2 \hat{\alpha}} \right) > 0, \label{eq:Jacdet} \noeqref{eq:Jacdet} \\[1ex]
\Delta(J) &= (\hat{r}^*_Y - \hat{\nu}_1 \hat{r}^*_X)^2 + \frac{\hat{\nu}_1 \hat{\gamma}^2}{\hat{\alpha}}\left( \frac{\hat{\nu}_1 \hat{\gamma}^2}{4 \hat{\alpha}} - (\hat{r}^*_Y + \hat{\nu}_1 \hat{r}^*_X) \right) > 0. \label{eq:Jacdis} \noeqref{eq:Jacdis}
\end{align}
Since $\det > 0$, $\tr > 0$, and $\Delta > 0$, the equilibrium point is an unstable node that persists in the non-singular limit $\hat{\sigma}_6 > 0$. At the singular limit $\hat{\sigma}_6 = 0$, the equilibrium point $(0, Y_{eq})$ loses hyperbolicity.

\section{Preservation of key geometric features}
\label{app:sec:geom}

\subsection{Compact domain and confined fold points}
\label{sec:theoderiv}

The relevant biophysical region of the $Y$-nullcline requires a non-negative denominator:
\begin{equation}
X^2 = \frac{\hat{\gamma} Y^\frac{1}{2} \left(\hat{\sigma}_3 Y + \hat{\sigma}_4 \right)}{-\hat{\gamma} \hat{\sigma}_1 Y^\frac{3}{2} + Y - \hat{\gamma} \hat{\sigma}_2 Y^\frac{1}{2} + \hat{\sigma}_6}\,.
\label{eq:ynullmnfd} \noeqref{eq:ynullmnfd}
\end{equation}
The zeros (poles) of the denominator are given by:
\begin{equation}
Y^\frac{3}{2} - \frac{1}{\hat\gamma \hat\sigma_1} Y + \frac{\hat{\sigma}_2}{\hat\sigma_1} Y^\frac{1}{2} - \frac{\hat{\sigma}_6}{\hat\gamma \hat\sigma_1} = 0\,.
\label{eq:ynullpols} \noeqref{eq:ynullpols}
\end{equation}
In the asymptotic limit $\hat\sigma_6 \to 0$, dominant balance for the root near $Y=0$ yields:
\begin{equation}
Y \approx \left( \frac{\hat\sigma_6}{\hat\gamma \hat\sigma_2} \right)^2\,.
\label{eq:lfold} \noeqref{eq:lfold}
\end{equation}
For this root to be asymptotically small, we require the scaling relationship:
\begin{equation}
\hat\sigma_2 \sim \mathcal{O}(\hat\sigma_6^a), \quad \text{for } 0 < a < 1.
\label{eq:rels6s2} \noeqref{eq:rels6s2}
\end{equation}
The maximum $Y$-value is bounded by the upper root, $Y_{max} \approx (1/(\hat\gamma \hat\sigma_1))^2 - 2(\hat\sigma_2/\hat\sigma_1)$. 

The S-shaped $Y$-nullcline has two fold points defined by the quartic:
\begin{equation}
Y^2 + 2 \hat\gamma \left(\frac{\hat\sigma_1 \hat\sigma_4}{\hat\sigma_3} - \hat\sigma_2 \right) Y^\frac{3}{2} + \left(3 \hat\sigma_6 - \frac{\hat\sigma_4}{\hat\sigma_3} \right) Y + \frac{\hat\sigma_4 \hat\sigma_6}{\hat\sigma_3} = 0\,.
\label{eq:ynullfdqr} \noeqref{eq:ynullfdqr}
\end{equation}
Assuming the small parameters scale equally:
\begin{equation}
\frac{\hat\sigma_4}{\hat\sigma_3} \sim \mathcal{O}(1)\,,
\label{eq:rels4s3} \noeqref{eq:rels4s3}
\end{equation}
the quartic has a non-trivial approximation for the upper fold ($Y_{f2}$):
\begin{equation}
Y_{f2} + 2 \hat\gamma \frac{\hat\sigma_1 \hat\sigma_4}{\hat\sigma_3} Y_{f2}^\frac{1}{2} - \frac{\hat\sigma_4}{\hat\sigma_3} \approx 0\,.
\label{eq:yf2as} \noeqref{eq:yf2as}
\end{equation}
This yields a unique positive root $Y_{f2} \sim \mathcal{O}(1)$:
\begin{equation}
Y_{f2}^\frac{1}{2} =  \sqrt{\left( \frac{\hat\gamma \hat{\sigma}_1 \hat{\sigma}_4}{\hat{\sigma}_3} \right)^2 + \frac{\hat{\sigma}_4}{\hat{\sigma}_3}} -  \frac{\hat\gamma \hat\sigma_1 \hat{\sigma}_4}{\hat{\sigma}_3}\,.
\label{Yf2-root}
\end{equation}
Substituting \eqref{Yf2-root} into \eqref{eq:ynullmnfd} shows $X_{f2} \sim \mathcal{O}(\hat\sigma_3^{1/2}) \to 0$.

For the lower fold, the lowest-order terms of \eqref{eq:ynullfdqr} provide the approximation:
\begin{equation}
Y_{f1} \approx \hat\sigma_6\,.
\label{Yf1-root}
\end{equation}
Substituting \eqref{Yf1-root} into \eqref{eq:ynullmnfd} gives:
\begin{equation}
X_{f1}^2 \approx \frac{\hat\gamma \hat\sigma_4 \sqrt{\hat\sigma_6}}{2 \hat\sigma_6 - \hat\gamma \hat\sigma_2 \sqrt{\hat\sigma_6}}\,.
\label{eq:xf1}
\end{equation}
For $X_{f1}$ to be non-negative and bounded away from zero, the numerator and denominator must balance. This refines \eqref{eq:rels6s2} to:
\begin{equation}
\hat\sigma_2 \sim \mathcal{O}(\hat\sigma_6^a): \frac{1}{2} \le a < 1\,,
\label{eq:s2minexp} \noeqref{eq:s2minexp}
\end{equation}
and requires the asymptotic scaling:
\begin{equation}
\hat\sigma_4 \sim \mathcal{O}(\hat\sigma_6^\frac{1}{2}) \Rightarrow \hat\sigma_3 \sim \mathcal{O}(\hat\sigma_6^\frac{1}{2})\,.
\label{eq:reqs3s4} \noeqref{eq:reqs3s4}
\end{equation}

\subsection{The unstable node}
\label{sec:unsnode}

The $X$-nullcline is defined by:
\begin{equation}
X^2 = \frac{\hat\alpha (\hat\sigma_3 Y + \hat\sigma_4)}{(\hat\nu_1 - \hat\alpha \hat\sigma_1)Y + (\hat\nu_1 \hat\sigma_6 - \hat\alpha \hat\sigma_2)}\,.
\label{eq:xnullmfd}
\end{equation}
Equating \eqref{eq:ynullmnfd} and \eqref{eq:xnullmfd} yields the unique equilibrium:
\begin{equation}
(X_{eq},Y_{eq})= \left( \left( \dfrac{\hat{\alpha} \left( \hat{\alpha}^2 \hat{\sigma}_3 + \hat{\gamma}^2 \hat{\nu}_1^2 \hat{\sigma}_4 \right)}{\hat{\gamma}^2 \hat{\nu}_1^2 (\hat{\nu}_1 \hat{\sigma}_6 - \hat{\alpha} \hat{\sigma}_2) + \hat{\alpha}^2 (\hat{\nu}_1 - \hat{\alpha} \hat{\sigma}_1 )} \right)^{1/2},  
\left(\dfrac{\hat{\alpha}}{\hat{\gamma} \hat{\nu}_1}\right)^2 \right)\,.
\label{eq:eqbmxy} \noeqref{eq:eqbmxy}
\end{equation}
For this to be an unstable node, it must lie on the middle branch ($Y_{f1} < Y_{eq} < Y_{f2}$). For $Y_{eq}$ to be bounded away from zero, we require:
\begin{equation}
\frac{\hat\nu_1}{\hat\alpha} \sim \mathcal{O}(1)\,.
\label{eq:reqalnu} \noeqref{eq:reqalnu}
\end{equation}
Substituting \eqref{eq:eqbmxy} and \eqref{Yf2-root} into the condition $Y_{eq} < Y_{f2}$ yields:
\begin{equation}
\left(\frac{\hat\alpha}{\hat\gamma\hat\nu_1}\right)^2 \left(1-2\hat\sigma_1\frac{\hat\alpha}{\hat\nu_1} \right)^{-1} < \frac{\hat\sigma_4}{\hat\sigma_3}\,,
\label{cond:param1}
\end{equation}
which necessitates the final condition (satisfied by \cref{tab:sumparams}):
\begin{equation}
\frac{\hat{\nu}_1}{\hat{\alpha}}> 2\hat{\sigma}_1\,.
\label{cond:param2}
\end{equation}

\section{Blow-up analysis}
\label{app:blowup}

This appendix collects the desingularised systems, linearisations, and the trapping-wedge estimate underlying \cref{sec:blowup}. 

\subsection{First cylindrical blow-up (Chart $K_1$)}
\label{app:bu1}\label{app:k1}

\paragraph{Invariant plane $\varepsilon_1 = 0$.} The system is:
\begin{equation}
r_1' = 0, \qquad Z_1' = (1 - \gamma \sigma_1 Z_1) Z_1^2.
\label{eq:sys:e1cha1}
\end{equation}
Stability is determined by $\lambda(Z_1) = (2-3\gamma \sigma_1 Z_1)Z_1$. $\Gamma_0^D$ is degenerate ($\lambda=0$), while $\Gamma_2$ is normally hyperbolic and attracting ($\lambda = -1/(\gamma \sigma_1) < 0$).

\paragraph{Invariant plane $r_1 = 0$.} The system is:
\begin{equation}
\begin{pmatrix} Z_1' \\ \varepsilon_1' \end{pmatrix} = \begin{pmatrix} (1 - \gamma \sigma_1 Z_1) Z_1 - \varepsilon_1^2 \gamma (\sigma_3 Z_1^2 + 1) \\ 2 \varepsilon_1^2 \left[ (\nu - \alpha \sigma_1) Z_1^2 - \varepsilon_1^2 \alpha (\sigma_3 Z_1^2 + 1) \right] \end{pmatrix} Z_1.
\label{eq:sys:r1cha1}
\end{equation}
At $p_1$, there is a stable eigenvalue $-1/(\gamma\sigma_1)$ (along $Z_1$) and a zero eigenvalue (along $\varepsilon_1$). The local centre manifold $W^c_{loc}(p_1)$ and its reduced flow are:
\begin{align}
Z_1 &= \frac{1}{\gamma \sigma_1} -  \frac{\sigma_3 + \gamma^2 \sigma_1^2}{\gamma \sigma_1^2} \varepsilon_1^2 + \mathcal{O}(\varepsilon_1^3), \label{app:eq:cmgraph} \\
\varepsilon_1' &= \frac{2(\nu - \alpha \sigma_1)}{\gamma^3 \sigma_1^3} \varepsilon_1^2 + \mathcal{O}(\varepsilon_1^3). \label{app:eq:cmflow}
\end{align}

\paragraph{The trapping wedge.} $W^c_{loc}(p_1)$ is geometrically trapped in a narrowing wedge in the $(\varepsilon_1, Z_1)$-plane. It lies strictly below the tangent line $L := \{Z_1 = 1/(\gamma\sigma_1)\}$ due to its quadratic departure, which is consistent with the vector field on $L$ \eqref{eq:sys:r1cha1} forcing trajectories into the lower half-plane ($Z_1' < 0$). Furthermore, it remains strictly above the $Z_1$-nullcline: the tangency condition $Z_1' = (dZ_1/d\varepsilon_1)\varepsilon_1'$ yields $Z_{1,\text{cm}} - Z_{1,\text{null}} > 0$ to leading order $\mathcal{O}(\varepsilon_1^3)$. Pinched between $L$ and the nullcline, the flow transitions towards the interior ($Z_1$ decreasing).

\subsection{Second cylindrical blow-up (Charts $K_3$ and $K_4$)}
\label{app:bu2}\label{app:k3}\label{app:k4}

\paragraph{Chart $K_3$ ($s_3 = 0$).} The invariant plane flow is:
\begin{equation}
Z_3' = Z_3^2 - \gamma (\sigma_2 r_3^2 + 1) Z_3 + r_3^4.
\label{app:eq:z3quad}
\end{equation}
At $p_4$, normal hyperbolicity is lost as the non-trivial eigenvalue $\lambda = 2 Z_3 - \gamma (\sigma_2 r_3^2+1)$ vanishes. 

\paragraph{Chart $K_3$ ($r_3 = 0$).} The invariant plane flow is:
\begin{equation}
\begin{pmatrix} Z_3' \\ s_3' \end{pmatrix} = \begin{pmatrix} Z_3^2 - \gamma Z_3 - s_3^2 \gamma Z_3^3 (\sigma_1 + s_3^2 \sigma_3) + 4 \alpha s_3^3 Z_3^2 + \dots \\ -2 \alpha Z_3 s_3^4 + \mathcal{O}(s_3^6) \end{pmatrix}.
\label{app:eq:r3cha3}
\end{equation}
Linearisation at $p_7$ yields $\lambda_1 = \gamma > 0$ and $\lambda_2 = 0$ (with centre flow $s_3' \approx -2 \alpha \gamma s_3^4$). The Jacobian at $p_3$ is $\text{diag}(0, -\gamma, 0)$.

\paragraph{Chart $K_4$ ($r_4 = 0$).} The invariant plane flow is:
\begin{equation}
\begin{pmatrix} s_4' \\ \varepsilon_4' \end{pmatrix} = \begin{pmatrix} s_4 \left( 1 - \gamma \sigma_1 s_4^2 - \varepsilon_4^2 \gamma (\sigma_3 s_4^4 + 1) \right) \\ -\varepsilon_4 \left( 1 - \gamma \sigma_1 s_4^2 - \varepsilon_4^2 \gamma (\sigma_3 s_4^4 + 1) - 4 s_4^5 \varepsilon_4 (\nu - \alpha \sigma_1) + \dots \right) \end{pmatrix}.
\label{app:eq:r4cha4}
\end{equation}
The Jacobian at the origin $p_8$ is $J(p_8) = \text{diag}(1, -1)$, confirming it is a hyperbolic saddle.

\paragraph{Chart $K_4$ ($s_4 = 0$).} The invariant vertical dynamics and the folded critical manifold $\Gamma_4$ are:
\begin{align}
\varepsilon_4' &= - \varepsilon_4 [1 - \varepsilon_4^2 \gamma (\sigma_2 r_4^2 + 1) + \varepsilon_4^4 r_4^4], \label{app:eq:s4cha4} \\
0 &= r_4^4 (\varepsilon_4^2)^2 - \gamma (\sigma_2 r_4^2 + 1) \varepsilon_4^2 + 1. \label{app:eq:g4quart}
\end{align}

\section{Geometric Singular Perturbation Theory (GSPT)}
\label{app:gspt}

We briefly describe the coordinate-independent setup of GSPT \cite{Wechselberger2020}. Consider the singular perturbation problem
\begin{equation}
x' = F(x; \varepsilon) = F_0(x) + \sum_{i \in \mathbb{N}} \varepsilon^i F_i(x), \quad x \in U_0 \subset \mathbb{R}^n, \quad n \geq 2\,,
\label{eq:app:gspt}
\end{equation}
where $t$ is the fast time ($'=d/dt$). The leading-order fast system is the \emph{layer problem} $x' = F_0(x)$ \label{eq:app:layer}\label{df:app:spt}\label{as:app:s0}\label{df:app:gslp}. The \emph{critical manifold} $S$ is a $k$-dimensional regular level set of equilibria, factored as $F_0(x) = N_0(x) f_0(x)$ \label{eq:app:ffac}\label{as:app:facnf}, where $f_0(x) = \mathbb{O}_{n-k}$. This induces a bundle splitting $T_x \mathbb{R}^n = T_x S \oplus \mathcal{N}_x$ \label{eq:app:gspl}, with $\mathcal{N}_x = \text{im}(N_0)$ representing the linear fast fiber subspace. $S$ is \textit{normally hyperbolic} if the matrix $Df_0 N_0|_S$ is hyperbolic \label{df:app:gsnh}.

Rescaling to the slow time $\tau = \varepsilon t$ yields the slow system $\dot{x} = \varepsilon^{-1} N_0 f_0 + F_1 + \mathcal{O}(\varepsilon)$ \label{eq:app:gssl}. The \emph{reduced problem} on $S$ is given by the oblique projection onto the tangent space $T_x S$ along $\mathcal{N}_x$:
\begin{equation}
\dot{x} = \Pi^S(x) F_1 (x),
\label{eq:app:gsrp}
\end{equation}
\label{rmk:app:oliq}
where the projection operator is explicitly defined as:
\begin{equation}
\Pi^S(x) = \mathbb{I}_n - \Pi^N(x) = \mathbb{I}_n - N_0(x) \left( Df_0(x) N_0(x) \right)^{-1} Df_0(x).
\label{eq:app:gspo}
\end{equation}

\begin{figure}[h]
\centering
\includegraphics[width=0.65\textwidth]{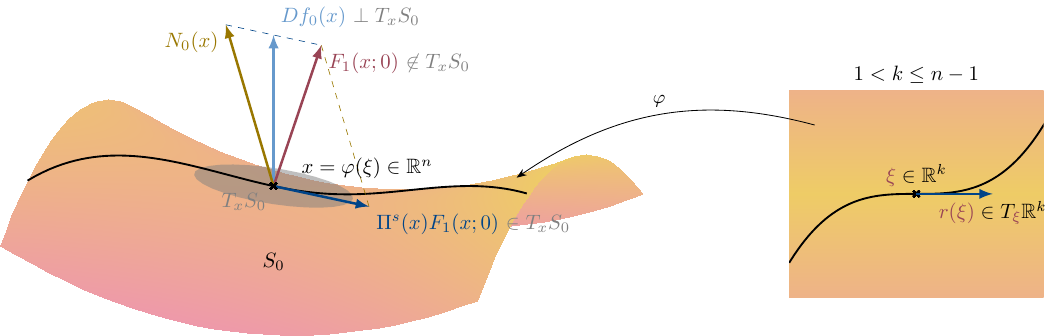}
\caption{\textit{Left:} Geometric representation of the reduced vector field obtained via the oblique projection $\Pi^S$ onto the tangent bundle $TS$ along the linear fast fiber bundle $\mathcal{N}$. \textit{Right:} Smooth embedding $\varphi$ from $U_1$ into $U_0$.}
\label{fig:app:gsptpm}
\end{figure}

Fenichel theory breaks down at \textit{contact points} $G$ where $\rk(Df_0 N_0) = n-k-1$ \label{df:app:Gnh}\label{lm:app:dsg}. The reduced flow has a finite-time blow-up here, but its direction can be evaluated using the smooth extension provided by the adjugate operator $\adj(A) = \det(A) A^{-1}$ \label{eq:app:adj_inv}\label{df:app:adjugate}\label{rmk:app:flch}. A contact point is a \textit{regular jump point} if the desingularised flow is non-zero: $N_0(x) \adj \left( Df_0(x) N_0(x) \right) Df_0(x) F_1(x) \neq \mathbb{O}_n$ \label{df:app:rjp}. At such points (e.g., the fold $p_4$), trajectories are ejected into the fast layer flow.

\section{Parametrisation method}
\label{app:pm}

To evaluate the reduced vector fields on higher-order slow manifolds, we use the \emph{parametrisation method} \cite{Lizarraga2021}. Adopting the setup from \cref{app:gspt}, we assume the $k$-dimensional normally hyperbolic critical manifold $S$ is the image of a smooth embedding \label{as:app:pm_nh}
\begin{equation}
\varphi_0 \colon U_1 \subset \mathbb{R}^k \hookrightarrow U_0 \subset \mathbb{R}^n\,.
\label{eq:app:emb0} \noeqref{eq:app:emb0}
\end{equation}

To find the slow flow on the perturbed manifold $S^\varepsilon$, we seek a smooth embedding $\varphi(\xi; \varepsilon) = \sum \varepsilon^i \varphi_i(\xi)$ and a reduced vector field $r(\xi; \varepsilon) = \sum \varepsilon^i r_i(\xi)$ such that $x(t) := \varphi(\xi(t); \varepsilon)$ satisfies $x' = F(x; \varepsilon)$ whenever $\xi' = r(\xi; \varepsilon)$. \label{app:pm:emred} This requires solving the \textit{conjugacy equation}:
\begin{equation}
\mathcal{F}(\varphi, {r})(\xi, \varepsilon) = D_\xi \varphi (\xi, \varepsilon) \circ {r}(\xi; \varepsilon) - {F}(\varphi (\xi, \varepsilon), \varepsilon) = 0.
\label{eq:app:pmcnj} \noeqref{eq:app:pmcnj}
\end{equation}

Matching powers of $\varepsilon$ in the Taylor expansion of \eqref{eq:app:pmcnj} yields a hierarchy of infinitesimal conjugacy equations. Utilising the factorisation $F_0 = N_0 f_0$, the equation for the $j$-th order pair $(\varphi_j, r_j)$ is:
\begin{equation}
D_\xi \varphi_0 \circ {r}_j - N_0(\varphi_0) D_x f_0(\varphi_0) \circ \varphi_j = {G}_j,
\label{eq:app:infconj}
\end{equation}
where $G_j$ depends only on known lower-order terms $\{ \varphi_i, r_i \}_{i < j}$. 

Applying the oblique projections $\Pi^S$ and $\Pi^N$ isolates the tangential and normal components, respectively. Since $D_\xi \varphi_0 \circ r_j \in T_x S$ and $N_0 D f_0 \circ \varphi_j \in \mathcal{N}_x$, we obtain the explicit update formulas for the internal dynamics and the slow manifold correction:
\begin{align}
r_j(\xi) &= (D_\xi \varphi_0)^L \left[ \Pi^S(\varphi_0(\xi)) \, G_j(\xi) \right], \label{eq:app:rj_sol} \\
\varphi_j(\xi) &= - (D_x f_0)^R \left( D_x f_0 N_0 \right)^{-1} D_x f_0 \, G_j(\xi). \label{eq:app:phij_sol}
\end{align}
Here, $(D_\xi \varphi_0)^L$ and $(D_x f_0)^R$ denote the left and right inverses of $D_\xi \varphi_0$ and $D_x f_0$. The specific choice of these generalised inverses determines structural properties of the manifold (e.g., maintaining a coordinate graph parametrisation versus shifting purely along fast fibers). Any tangential component introduced by these choices represents a reparametrisation of the manifold, which is automatically compensated by a corresponding shift in the reduced vector field $r_j$ to satisfy \eqref{eq:app:infconj}.

\section{Persistence of the Singular Limit Cycle}
\label{sec:returnmap}

To prove Proposition \ref{prop:persistence}, we construct a return map $P_\varepsilon$ on a cross-section $\Sigma$ in the leaf $L_\varepsilon$, and bound its derivative $\partial_h P_\varepsilon$. We partition the singular cycle $\Gamma^0$ into seven local passages. By the chain rule, the overall transverse contraction is bounded by the product of the local bounds: 
\[
|\partial_h P_\varepsilon| \le \prod_{j=1}^7 B_j.
\]
The seven passages fall into four categories of established local dynamics:

\paragraph{Normally hyperbolic slow segments and the fold.} 
The passages along the three slow segments ($\Gamma_2, \Gamma_5, \Gamma_{4a}$) are governed by standard Fenichel theory \cite{Fenichel1979}. Since the transverse eigenvalues are negative and $\mathcal{O}(1)$, the transit time $\mathcal{O}(1/\varepsilon)$ gives an exponential contraction bound $B_j = \mathcal{O}(e^{-c_j/\varepsilon})$ for $j=1,3,5$ with $c_j>0$. These passages remain within a Hausdorff distance of $\mathcal{O}(\varepsilon)$ from $\Gamma^0$. 
Similarly, the passage past the generic fold $p_4$ is governed by standard planar fold theory \cite{Krupa2001}, which contributes a further exponential contraction factor $B_6 = \mathcal{O}(e^{-c_6/\varepsilon})$ and a maximal deviation of $\mathcal{O}(\varepsilon^{2/3})$ from the singular limit.

\paragraph{Centre-manifold gateways.} 
At the gateways $p_1$ and $p_3$, normal hyperbolicity is lost. However, the spectra contain exactly one negative eigenvalue and zeros, with no unstable directions. The local centre manifolds are therefore normally attracting. Thus, transverse distances are non-expanding, contributing $\mathcal{O}(1)$ bounded factors $B_2, B_4 \le K$, and the trajectories remain within $\mathcal{O}(\varepsilon)$ of $\Gamma^0$.

\paragraph{The equatorial saddle.} 
At $p_5$, the blow-up geometry yields an equatorial saddle. In the transverse plane, the point acts as a hyperbolic saddle with hyperbolicity ratio $\rho=1$. Applying the standard resonant Dulac map, the transition between transversals yields at most a logarithmic stretching factor $B_7 = \mathcal{O}(\ln(1/\varepsilon))$. 
Importantly, because $p_5$ lies on the equator of the cylinder rather than at a singular corner of the blow-up, the parameter $\varepsilon$ is tied directly to the standard radial coordinate. Consequently, the trajectory remains within $\mathcal{O}(\varepsilon)$ of the singular limit, avoiding the fractional distance penalty associated with corner passages.

\paragraph{Contraction and convergence rate.}
Substituting the local bounds into the chain rule product yields:
\[
|\partial_h P_\varepsilon| \le K^2 \cdot \mathcal{O}(\ln(1/\varepsilon)) \cdot e^{-c/\varepsilon}, \qquad c = c_1 + c_3 + c_5 + c_6 > 0.
\]
Since the exponential contraction $e^{-c/\varepsilon}$ from the slow manifolds and fold dominates the logarithmic expansion $\mathcal{O}(\ln(1/\varepsilon))$ from the saddle, for any $0 < C < c$ there exists $\varepsilon_0 > 0$ such that $|\partial_h P_\varepsilon| \le e^{-C/\varepsilon} < 1$ for all $\varepsilon \in (0, \varepsilon_0]$. $P_\varepsilon$ is therefore a strict contraction, guaranteeing the existence, uniqueness, and exponential stability of the limit cycle $\Gamma_\varepsilon$.
Finally, the overall Hausdorff distance $\dist_{\mathrm H}(\Gamma_\varepsilon, \Gamma^0)$ is bounded by the maximum of the local deviations. Comparing $\mathcal{O}(\varepsilon)$ from the slow segments, gateways, and saddle against $\mathcal{O}(\varepsilon^{2/3})$ from the fold, the fold passage forms the bottleneck, establishing the global convergence rate of $\mathcal{O}(\varepsilon^{2/3})$ and completing the proof.

\end{appendices}

\bibliography{MSPaper_abbrev}

\end{document}